\documentclass[11pt]{article}
\usepackage{amssymb,amsmath,amsthm,lmodern,cancel,graphicx,calrsfs}
\usepackage[T1]{fontenc}
\usepackage[latin9]{inputenc}
\usepackage{textcomp}
\usepackage[dvipsnames]{xcolor}
\usepackage{enumerate}
\usepackage{thmtools, mathtools} 
\usepackage[colorlinks=true,hyperindex=true]{hyperref}

\def\titlerunning#1{\gdef\titrun{#1}}
\makeatletter
\def\author#1{\gdef\autrun{\def\and{\unskip, }#1}\gdef\@author{#1}}
\def\address#1{{\def\and{\\\hspace*{18pt}}\renewcommand{\thefootnote}{}%
\footnote {#1}}%
\markboth{\autrun}{\titrun}}
\makeatother
\def\email#1{e-mail: #1}
\def\subjclass#1{{\renewcommand{\thefootnote}{}%
\footnote{\emph{Mathematics Subject Classification (2010):} #1}}}
\def\keywords#1{\par\medskip
\noindent\textbf{Keywords.} #1}

\usepackage[activate={true,nocompatibility},final,tracking=true,kerning=true,spacing=true,factor=1100,stretch=15,shrink=15]{microtype}
\microtypecontext{spacing=nonfrench}

\newtheorem{thm}{Theorem}[section]
\newtheorem{cor}[thm]{Corollary}
\newtheorem{lem}[thm]{Lemma}
\newtheorem{prop}[thm]{Proposition}

\newtheorem{mainthm}[thm]{Main Theorem}

\newtheorem{theorem}{Theorem}

\newtheorem{lemma}[theorem]{Lemma}

\theoremstyle{definition}

\numberwithin{equation}{section}

\newcommand{\abs}[1]{\left\lvert #1 \right \rvert}

\newcommand{\norm}[1]{\lVert #1 \rVert}

\newcommand{\mc}[1]{\mathcal{#1}}
\newcommand{\m}[1]{\mathbb{#1}}

\newcommand{\rar}[0]{\rightarrow}

\def\st{, \,}
\def\ie{i.e. }

\renewcommand\Re{\operatorname{Re}}
\renewcommand\Im{\operatorname{Im}}

\def\a{\alpha}
\def\b{\beta}
\def\g{\gamma}

\def\d{\delta}
\def\D{\Delta}

\def\t{\theta}

\def\k{\kappa}
\def\l{\lambda}

\def\s{\sigma}

\def\vare{\varepsilon}
\def\HH{{\mathbb H}}
\def\eps{\varepsilon}
\def\Lip{{\rm Lip}}
\def\Rev{{\rm Rev}}

\def\dd{\,\mathrm{d}}

\newcommand{\ad}[1]{\underline{#1}}

\def\P{\mathrm P}

\newcommand\SLE{\operatorname{SLE}}

\begin{document}




\titlerunning{Reversibility of the Loewner energy via SLE$_{0+}$}

\title{The energy of a deterministic Loewner chain: Reversibility and interpretation via SLE$_{0+}$}

\author{Yilin Wang}

\date{March 5, 2019}

\maketitle

\address{ETHZ, D-Math, Raemistrasse 101, 8092 Zurich, Switzerland; \email{yilin.wang@math.ethz.ch}}

\subjclass{Primary 30C55; Secondary 60J67}


\begin{abstract}
We study some features of the energy of a deterministic chordal Loewner chain, which is defined as the Dirichlet energy of its driving function. 
In particular, using an interpretation of this energy as a large deviation rate function for SLE$_\kappa$ as $\kappa \rar 0$ and the known reversibility of the SLE$_\kappa$ curves for small $\kappa$, we show that the energy of a deterministic curve from one boundary point $a$ of a simply connected domain $D$ to another boundary point $b$, 
 is equal to the energy of its time-reversal ie. of the same curve but viewed as going from $b$ to $a$ in $D$. 

\keywords{Loewner differential equation, Loewner energy, reversibility, quasiconformal mapping, Schramm-Loewner Evolution.}
\end{abstract}

\section{Introduction}

The Loewner transform is a natural way to encode certain two-dimensional paths by a one-dimensional continuous function, via a procedure that involves 
iterations of conformal maps. This idea has been introduced by Loewner in 1923 \cite{Loewner1923} in its radial form (motivated by the Bieberbach conjecture, and it eventually was an important tool in its solution \cite{DeBranges1985}, \cite{hayman1994multivalent}). It has also a natural chordal counterpart, when one encodes a simple curve from one boundary point to another instead of a simple curve from a boundary point to some inside point as in the radial case.

Let us briefly recall this chordal Loewner description of a continuous simple curve $\gamma$ from $0$ to infinity in the open upper half-plane $\HH$. We parameterize the curve in such a way that 
the conformal map $g_t$ from $\HH \setminus \gamma [0,t]$ onto $\HH$ with $g_t (z)=  z + o(1)$ as $z \to \infty$ satisfies in fact $g_t ( z)= z + 2t / z + o( 1/z)$ (it is easy to see that it is always possible to reparameterize a continuous curve in such a way). It is also easy to check that one can extend $g_t$ by continuity to the  boundary point $\gamma_t$ and that the real-valued function $W_t := g_t ( \gamma_t)$ is continuous. Loewner's equation then shows that the function $t \mapsto W_t$ (that is often referred to as the driving function of $\gamma$) does fully characterize the curve. 

The Loewner energy (we will also often just say ``the energy'') of this Loewner chain is then defined to be the usual Dirichlet energy of its driving function (this quantity has been also introduced and considered in the very recent paper by Friz and Shekhar \cite {friz2015}): 
$$ I(\gamma) := \frac 1 2 \int_0^{\infty} \left( \frac {dW_t}{dt} \right)^2 dt = \frac 1 2 \int_0^{\infty} \left( \frac {d ( g_t ( \gamma_t))}{dt} \right)^2 dt .$$ 
This quantity does not need to be finite, but of course, for a certain class of sufficiently regular curves $\gamma$, $I(\gamma)$ is finite. 
Conversely, as we shall explain, any real-valued function $W$ with finite energy does necessarily correspond to a simple curve $\gamma$. 

It is immediate to check that if one considers the simple curve $u \gamma$ instead of $\gamma$ for some $u > 0$, then $I (u \gamma) = I ( \gamma)$. This scale-invariance property of the energy 
makes it possible 
to define the energy $I_{D,a,b} (\eta)$ of any simple curve $\eta$  from a boundary point $a$ of a simply connected domain $D$ to another boundary point $b$ to be the energy of the conformal 
image of $\eta$ via any uniformizing map $\psi$ from $(D,a,b)$ to $(\HH, 0, \infty)$ (as it will not depend on the actual choice of $\psi$). 

For such a simple curve $\eta$, one can define its time-reversal $\hat \eta$ that has the same trace as $\eta$, but which is viewed as going from $b$ to $a$. 
The first main contribution of the present paper is to prove the following reversibility result: 
\begin {mainthm}
\label {main1}
The Loewner energy of the time-reversal $\hat \eta$ of a simple curve $\eta$ from $a$ to $b$ in $D$ is equal to the Loewner energy of $\eta$: $I_{D,a,b}(\eta) = I_{D,b,a}(\hat \eta)$. 
\end {mainthm}

Note that if $\gamma$ is a curve from $0$ to infinity in the upper half-plane, then $-1/\gamma$ is a curve from $\infty$ to $0$ in the upper half-plane, but if we view it after time-reversal, one gets
again curve from $0$ to infinity.  An equivalent formulation of the theorem is that the Dirichlet energy of the driving function of this curve is necessarily equal to that of the driving function of $\gamma$.

A naive guess would indicate that this simple statement must have a straightforward proof, but this seems not to be the case.
Indeed, the definition of the curve $\gamma$ out of $W$ via iterations of 
conformal maps is very ``directional'' and not well-suited to time-reversal. As it turns out, the way in which the energy of the driving function $W$ of $\gamma$ will be redistributed on the driving function of its time-reversal is highly non-local and rather intricate. This suggest that a more symmetric reinterpretation of the energy of the Loewner chain will be handy in order to prove this theorem -- and this is indeed the strategy that we will follow. 

The reader acquainted to Schramm Loewner Evolutions (SLE) may have noted features reminiscent of part of the SLE theory. Recall that chordal SLE$_\kappa$ as defined by Schramm in \cite{schramm2000scaling} is
the random curve that one obtains when one chooses the driving function to be $\sqrt {\kappa}$ times a one-dimensional Brownian motion. This is of course a function with infinite energy, but 
 Brownian motion and Dirichlet energy are not unrelated  (and we will comment on this a few lines).
Recall also from \cite{Schramm2005basic} that (and this is a non-trivial fact) chordal SLE$_\kappa$ is indeed a simple curve when $\kappa \in [0,  4]$ (see also Friz and Shekhar \cite {friz2015} for an approach based on rough paths).
Furthermore, given that  SLE$_\kappa$ curves should conjecturally arise as scaling limits of lattice based models from statistical physics, it was natural to conjecture that 
chordal SLE is reversible \ie that the time-reversal of a chordal SLE$_\kappa$ from $a$ to $b$ in $D$ is a chordal SLE$_\kappa$ from $b$ to $a$ in $D$ 
(note that here, this means an identity in distribution between two random curves). Proving this has turned out to be a challenge that resisted for some years, but has been settled by Zhan \cite{zhan2008} in the case of the simple curves $\kappa \in [0, 4]$, via rather non-trivial couplings of both ends of the 
path (see also Dub\'edat's commutation relations \cite{Dub{\'e}dat2007commutation}, Miller and Sheffield's approach based on the
Gaussian Free Field  \cite{MS2012imaginary1,MS2012imaginary2,MS2012imaginary3} that also provides a proof of 
reversibility for the non-simple case where $\kappa \in (4,8]$). 
In view of Theorem \ref {main1}, one may actually wonder if a possible roadmap towards proving reversibility of SLE is to deduce it from the reversibility of the Loewner energy for smooth curves $\gamma$. 
Our approach in the present paper will however be precisely the opposite; we will namely use the reversibility of the SLE$_\kappa$ curves as a tool to deduce the reversibility of the energy of smooth curves.

One way to relate Brownian motion to Dirichlet energy is to view the latter as the large deviation functional for paths with small Brownian component.
Recall that the law of one-dimensional Brownian motion on $[0,T]$ translated by a function $\l$ with $\l (0)=0$ is absolutely continuous with respect to the 
law of Brownian motion if and only if $ I_T (\l) := \int_0^T \dot{\l}(t)^2 \dd t / 2 $ is finite (the space of such function is called the Cameron-Martin space).  
In particular, if we set $I_T (\l) = \infty$ for all continuous functions in the complement of the Cameron-Martin space,
then $I_T$ is also the rate function for the large deviation principle on Brownian paths in the classical Schilder Theorem \cite{Zeitouni}. 
This loosely speaking means that $I_T (\l)$ measures the exponential decay rate for $\sqrt{\k}$ times a Brownian motion to be in a small neighborhood of the path $\l$. 

Our proof of Theorem \ref{main1} is based on the fact that the energy of the Loewner chain is the large deviation rate function of the SLE$_\kappa$ driving function when $\kappa \to 0+$. 
Loosely speaking, for a given $\gamma$ with finite energy, we want to relate $I( \gamma)$ to the decay as $\kappa \to 0$ of the probability  
for an SLE$_\kappa$ to be in a certain neighborhood $B_\vare (\gamma)$ of $\gamma$ via a formula of the type
$$I(\g) = \lim_{\vare \rar 0} \left(  \lim_{\k \rar 0} -\k \ln \P[ \rm{SLE}_\kappa \in  B_{\vare}(\g)]\right).$$
The main point will be to prove this for a reversible notion of $B_\vare$  
(so that a path is in the $\vare$-neighborhood of $\gamma$ if and only if its time-reversal is in the $\vare$-neighborhood of the time-reversal of $\gamma$). 
Once this will be done, the reversibility of the energy of the Loewner chain will follow immediately from this expression and the reversibility of SLE$_\kappa$ for all 
small $\kappa$. Similar idea is outlined by Julien Dub{\'e}dat in section~9.3 of \cite{Dub{\'e}dat2007commutation}.

Usual neighborhoods are not so well-suited for our purpose: Taking the neighborhood of $\g$ in the sense of some Hausdorff metric in the upper half-plane
is not well-adapted to the large deviation framework for the driving function, and the  $L^{\infty}$ norm on the driving function does not define reversible neighborhoods. 
Our choice will be to fix a finite number of points on the left side and the right side of $\g$, and to consider the collection of driving functions whose Loewner transform passes on the same side of these points than $\g$ does. This set of driving functions will play the role of $B_{\vare}(\g)$ and the large deviation principle will apply well (the $\vare \rar 0$ limit will be  replaced by letting the set of constraint points become dense in the upper half-plane). In order to properly apply the convergence, we will use some considerations about compactness of the set of $K$-quasiconformal curves, and the relation between finite energy chains and quasiconformal curves. 

In this paper, we will use basic concepts from the theory of quasiconformal maps, the Loewner equation, Schilder's large deviation theorem, SLE curves and some of their properties. We of course refer to the basic corresponding textbooks for background, but we choose to briefly and heuristically recall some of the basic definitions and results that we use, in order to help the readers in case some part of 
this background material does not belong to their everyday toolbox.
~\\

The paper is structured as follows: We will first derive in Section \ref {sec_finite_curve} some facts on deterministic Loewner chains with finite energy and their regularity. In Section \ref {sec_minimizing_curve_finite}, 
we derive Theorem \ref {main1} using the strategy that we have just outlined, and we conclude in Section \ref {commentssection} establishing some connections with ideas from SLE restriction properties and 
SLE commutation relations.  
In  upcoming  work, we plan to address generalizations and consequences of the present work to the case of full-plane curves, loops and more general graphs.

\section{Loewner chains of finite energy and quasiconformality} \label{sec_finite_curve}

\subsection{The Loewner transform}

 We say that  a subset $K$ of $\m{H}$ is a \emph{Compact $\m{H}$-Hull of half-plane capacity  $\rm{hcap}(K)$ seen from $\infty$}, if (i) $K$ is bounded, (ii) $H_K := \m{H} \backslash K$ is simply connected, and (iii) the unique conformal transformation $g_K : H_K \rar \m{H}$ such that $g_K(z) = z + o(1) $ when $z \to \infty$ also satisfies
$$g_K (z) = z + \frac{\rm{hcap}(K)}{z} + o(\frac{1}{z}),$$
when $z \to \infty$. We will refer to $g_K$ as the \emph{mapping-out function} of $K$.

Let $\mc{K}$ denote the set of compact $\m{H}$-hulls, endowed with the \emph{Carath\'eodory topology}, \ie a metric of the type
$$d_K (K_1, K_2)  = d(g_{K_1}^{-1}, g_{K_2}^{-1}),$$
where the metric $d$ generates the topology corresponding to  uniform convergence on compact subsets of the upper half-plane.
Recall  that this topology generated by $d_K$ differs from the Hausdorff distance topology. For instance, the arc $\{e^{i\t}, \t \in (0, \pi - \vare]\}$ converges to the half-disc of radius $1$ and center $0$ for the Caratheodory topology, but not for the Hausdorff metric. Let $(K_t)_{t>0} \subset \mc{K}$ be an increasing family of compact $\m{H}$-hulls for inclusion. 
For $s<t$, define $K_{s,t} = g_{K_s}(K_t \backslash K_s)$. We say that the sequence $(K_t)_{t>0}$ has \emph{local growth} if $\text{diam}(K_{t,t+h})$ converges to $0$ uniformly on compacts in $t$ as $h \to 0$
(where diam is the diameter for the Euclidean metric).

Let $\mc{L}_{\infty}$ denote the set of all increasing sequences of compact $\m{H}$-hulls $(K_t)_{t \in [0,\infty)}$ having local growth, parameterized in the way such that hcap$(K_t) =2t$. We endow $\mc{L}_{\infty}$ with the topology of uniform convergence on compact time-intervals.

Now we describe the chordal Loewner transform, which associates each continuous real-valued function with a family in $\mc{L}_{\infty}$: 
When $\l \in C([0,\infty))$, consider the Loewner differential equation 
\begin{equation} \label{eq_chordal_loewner}
 \partial_t g_t(z)  = 2/ (g_t(z) - \l_t)
\end{equation}
with initial condition $g_0(z)=z$.
The \emph{chordal Loewner chain in $\m{H}$ driven by the function} $\l$ (or \emph{the Loewner transform of $\l$}), is the increasing family $(K_t)_{t>0} $ defined by
\begin{align*}
K_t &= \{z \in \m{H} \st \tau(z) \leq t\},
\end{align*}
where $\tau(z)$ is the maximum survival time of the solution,
\ie 
$$\tau(z) = \sup\{ t \geq 0 \st \inf_{0\leq s \leq t} \abs{g_s(z)- \l_s} > 0\}.$$
It turns out that $(K_t)_{t>0}$ is indeed in $\mc{L}_{\infty}$. The family $(g_t)$ is sometimes called the \emph{Loewner flow} generates by $\l$.
We will also use $f_t := g_t - \l_t$, referred to as the \emph{centered Loewner flow}. 

Note finally that the Loewner equation is also defined for $z \in \m{R}$, and it is easy to see that the closure of $K_t$ in $\overline {\m{H}}$ satisfies $\overline {K_t} = \{z \in \overline{\m{H}} \st \tau(z) \leq t\}$.

The following result describes the explicit inverse of the Loewner transform and tells us that both procedures are continuous when we equip $C([0,\infty))$ with the topology of uniform convergence on compact intervals. (We will label all known results that we recall without proof 
by letters, and results derived in the present paper by numbers).

\begin{theorem}[see \cite{lawler2001values} 2.6]\label{thm_Loewner_Kufarev}
           The \emph{Loewner transform} $L: C([0,\infty)) \rar \mc{L}_{\infty}$ is a homeomorphism. The inverse transform is given by
$$\l_t = \bigcap_{h>0} \overline {K_{t,t+h}}, \quad \forall t \geq 0$$
where $K_{s,s'} = g_s(K_{s'} \backslash K_s)$ for $s<s'$.
\end{theorem}

We will also mention radial Loewner chains that are defined in a similar way in the unit disc $\m{D}$:
 The radial Loewner transform of $\xi \in C([0,\infty))$ is the $(K_t)_{t\geq 0}$ family given by the radial Loewner equation
$$ \partial_t h_t(z) = -h_t(z) \frac{h_t(z) +e^{i\xi(t)} }{ h_t(z) -e^{i\xi(t)} } $$ 
with initial condition $ h_0(z) =z$ for all $z$ in the closed unit disc,
where $K_t$ is defined as in the chordal case.

\subsection{The Loewner energy}\label{sec_Loewner_energy}

Let $T \in [0, \infty)$ and let $C_0([0,T])$ be the set of continuous functions $\l$ on $[0,T]$  with $\l(0) = 0$,  endowed with the $L^\infty$ norm. 
A function $\l\in C_0([0,T])$ is \emph{absolutely continuous}, if there exists a function $h \in L^1([0,T])$,  such that $\int_0^t h(s) ds = \l(t)$, $\forall t \geq 0$. In the sequel we write $\dot{\l} :=h$ and $\l$ is implicitly absolutely continuous whenever $\dot{\l}$ is considered. The \emph{energy up to time $T$} of the Loewner chain $\g$  driven by $\l$ is defined as 
$$I_T(\g) = I_T(\l) := 
\frac{1}{2}\int_0^T \dot{\l}(s)^2 \dd s$$ 
 for $\l$ absolutely continuous and $I_T ( \l) = \infty$ in all other cases. 
This energy has been recently introduced and used in the paper \cite {friz2015}, with the goal of 
providing a rough path approach to some features of SLE. While this SLE goal is quite different from the present paper, 
it has some similarities ``in spirit'' with the present paper, as it tries to 
provide a more canonical  (and less based on It\^o stochastic calculus -- therefore less directional) approach to SLE. Note that our project was developed independently of 
\cite {friz2015}. 

Let us list a few properties of this energy: 

\begin {itemize}
 \item 
The map $f \mapsto I_T (f)$ is lower-semicontinuous. 
Indeed, by elementary analysis, 
$I_T(\l)$
is the supremum over all finite partitions $\Pi = (0=t_0 < t_1<  t_2< \ldots< t_k = T)$ of
$$\sum_{t_{j-1},t_j \in \Pi} \frac{\abs{\l(t_j) - \l(t_{j-1})}^2}{t_j-t_{j-1}}.$$ 
Hence, if $\l_n$ is a sequence of functions converge uniformly on $[0,T]$ to $\lambda$, then
$$
\liminf_{n\rar \infty} I_T(\l_n)
\geq  \sup_{\Pi}  \ \liminf_{n\rar \infty} \sum_{t_{j-1},t_j \in \Pi} \frac{\abs{\l_n(t_j) - \l_n(t_{j-1})}^2}{t_j-t_{j-1}} \\
 =  I_T(\l).
$$

\item
The sets  
$\{\l \in  C_0([0,T]) \st I_T(\l) \leq c\}$
are compact for every $c \geq 0$ in $C_0([0,T])$.
This follows from the fact that a bounded energy set is equicontinuous and thus relatively compact in $C_0([0,T])$ by the Arzel\`a-Ascoli Theorem. 

\item
Similarly, we define \emph{the energy} of $\g$ driven by $\l$ as $I(\l):=I_{\infty}(\l)$ for $\l \in C_0 ([0, \infty))$ (endowed with the topology of uniform convergence on compact intervals). We still have that $I$ is lower semicontinuous, and the set of functions such that $I (\l ) \leq c$ is compact in $C_0 ([0, \infty))$.  Sometimes we will omit the subscript of $I_T$ if there is no ambiguity, and for $\l \in  C_0([0,T])$, we also define
$I(\l):=I_T(\l) = I(\l(\cdot \wedge T))$. 

\item We write $H_T \subset C_0([0,T])$ for the set of finite $I_T$ energy functions, and similarly, $H \subset C_0([0,\infty))$ for finite $I$ energy functions.

\item
For a simple curve $\gamma$ from $0$ to $\infty$ in the upper half-plane driven by the function $\l$, the energy satisfies obviously the following scaling property:  
$$I (\g) = I (u \g), \quad \forall u>0$$
because the driving function of $u\g$ is $\l^u (t) = u \l(t/u^2)$. 
Thus, the energy is invariant under conformal equivalence preserving $(\m{H}, 0, \infty)$. So one can define the Loewner energy of a curve from a boundary point $a$ to another boundary point $b$ in a simply connected domain $D$, by applying any conformal map from $(D,a,b)$ to $(\m{H}, 0 , \infty)$.  We use $I_{D,a,b}$ to indicate the energy in a different domain than $(\m{H}, 0 , \infty)$.

\item
It is also straightforward to see that the Loewner energy satisfies the simple additivity property
$$ I(\g) = I (f_s(\g[s, \infty))) + I_s(\g),$$
using the fact that the driving function of $f_t(\g[t, \infty))$ is $\l(s+t) - \l(s)$.

\item
If one is looking for a functional of a Loewner chain that satisfies additivity and conformal invariance, one is therefore looking for additive functionals of the driving functions with the right scaling property. There are of course 
several options (for instance one can add to $I_T$ the integral of the absolute value of the second derivative of $\lambda$ to the appropriate power). 
However, taking the integral of the square of the derivative is the most natural choice -- and it does satisfy the reversibility property that we derive in the present paper. 

\item 
Let us finally remark that a Loewner chain with $0$ energy is driven by the $0$ function, hence is the imaginary axis in $\m{H}$ which is the hyperbolic geodesic between $0$ and $\infty$. In view of the conformal invariance,  the energy of a curve from $a$ to $b$ in $D$ can therefore be viewed as a way to  measure how much does a curve differs from the hyperbolic geodesic from $a$ to $b$ in $D$ (see also related comments in Section~\ref{sec_conformal_restriction} for the conformal restriction properties). 

\end {itemize}

It may seem at first that the definition of the Dirichlet energy is quite ad hoc and depends on our choice of uniformizing domain for the Loewner flow (the fact that one works in the 
upper half-plane). Let us briefly indicate that this is not the case. Suppose for instance that we are working in the 
simply connected domain $D$, and that $\partial D$ is a smooth analytic arc in the neighborhood of $b$. It is then natural to parameterize the curve $\gamma$ by its capacity in $D$ as seen from $b$, \ie by 
$${\rm cap} (K):=  -\frac{1}{6} S \psi_K (b) = -\frac{1}{6} \left[\frac{\psi_K'''(b)}{\psi_K'(b)}- \frac{3}{2} \left(\frac{\psi_K''(b)}{\psi_K'(b)}\right)^2\right] = -\frac{1}{6} \left[\psi_K'''(b)- \frac{3}{2} \psi_K''(b)^2\right], $$
where $\psi_K : D \setminus K  \rar D$ is a conformal mapping such that $\psi_K(b)= b$ and $\psi_K' (b) = 1$. 
These derivatives are well defined due to the Schwarz reflection principle.
It is not hard to see that there exists a conformal equivalence $h$ from $(D,a,b)$ to $(\m{H},0,\infty)$ and $\alpha \in \m C^*$ such that $\alpha^2 {\rm cap} (\cdot)$ equals to the half-plane capacity of its image by $h$. 
For any continuous parametrization by $[0,T]$ of the oriented chord $\g$, let $\psi_t = \psi_{\g[0,t]}$ such that $\psi_t$ maps moreover $\g_t$ to $a$.  
Then the driving function  defined as $$\l(t) := \frac{ \a}{2} \psi_t''(b)  $$ corresponds to the driving function (parametrized by $\alpha^2 {\rm cap} (\g[0,t])$) of the image by $h$ of the curve in the half-plane.   
In other words, one can express the energy of the Loewner chain directly as the Dirichlet energy of $\psi_t'' (b)/2$ viewed as function of ${\rm cap} (\g [0,t])$ which are quantities naturally and directly defined in $D$.

\subsection{K-quasiconformal curve}

The following result will explain that a Loewner chain with finite energy is necessarily a simple curve: 
\begin{prop}\label{prop_finite_curve} For every $\l \in H$, there exists $K= K(I(\l))$, depending only on $I(\l)$, such that the trace of Loewner transform $\g$ of $\l$ is a $K$-quasiconformal curve, \ie the image of $i\m{R}_+$ by a $K$-quasiconformal mapping preserving $\m{H}$, $0$ and $\infty$.
\end{prop}
Properties of quasiconformal mappings are used in this section without proofs. For a panoramic survey on the topic, readers might refer to \cite{lehto2012univalent} Chapter~1, and detailed proofs 
can be found in \cite{lehto1973quasiconformal}.
Recall that a conformal mapping is locally a scalar times a rotation, which dilates a neighborhood of a point in all directions by a same constant.  A quasiconformal mapping is informally a orientation-preserving homeomorphism which locally does not stretch in one direction much more than in any other.   
One way to measure the deviation from being a conformal mapping, is to study how much the module $M(B)$ of a ring domain $B$ can change under the mapping (recall that it is preserved by conformal maps). Every ring domain $B$ is conformally equivalent to one of  the following annuli: $0 < \abs{z}< \infty$, $1 < \abs{z} < \infty $ or $1 < \abs{z}< R$. The \emph{module} $M(B)$ in the last case is defined as $\ln(R)$, and in the first two cases as infinite.
An orientation-preserving homeomorphism $f: D\subset{\m{R}^2} \rar E\subset{\m{R}^2}$ is \emph{$K$-quasiconformal}, if
$$M(f(B))/M(B) \in [1/K,K]   \text{ for every ring domain } B \text{ such that  } \overline {B} \subset D. $$ 
 In particular, conformal mappings are $1$-quasiconformal.

The following corollary of Proposition~\ref{prop_finite_curve} will be  useful in our proof of reversibility of the energy. 
A quasiconformal mapping $\varphi$ is said to be \emph{compatible} with $\l$ if it preserves $\m{H}, 0$ and $\infty$, and if $\varphi(i\m{R}_+) = \g$ and $\abs{\varphi (i)} = 1$.

\begin{cor}\label{cor_compact}
   Let $(\l_n)$ be a sequence of driving functions with energy bounded by $C<\infty$. Let $\varphi_n$ be a $K(C)$-quasiconformal mapping compatible with $\l_n$, there exists $\l \in H$ with $I(\l) \leq C$, and $\varphi$ compatible with $\l$, such that on a subsequence, $\l_{n}$ and $\varphi_{n}$ converge respectively (uniformly on compacts) to $\l$ and $\varphi$.  
\end{cor}

\begin{proof}
   Since the family of $K(C)$-quasiconformal mapping preserving $\m{H}$ is a normal family (cf.~\cite{lehto2012univalent} I.2.3),  together with the compactness of $I^{-1}([0,C])$, we can extract a subsequence where $\l_n$ and $\varphi_n$ converge both and respectively to $\l$ and $\varphi$. 
   The limit $\varphi$ is either a $K(C)$-quasiconformal map or a constant map into $\partial \m{H}$. The latter is excluded because of the choice $\abs{\varphi_n(i)} = 1$, and that the energy needed to touch a point near $\pm 1$ is not bounded (we will see it in Proposition~\ref{prop_gd}). Thus $\varphi$ is $K(C)$-quasiconformal, and we have to show that the quasiconformal curve $\tilde{\g}:=\varphi(i\m{R}_+)$ is the Loewner transform of $\l$.

   As explained in \cite{lind2010collisions}, Theorem 4.1. and Lemma 4.2, the modulus of continuity of $\g_n$ depends only on $K(C)$ when parameterized by capacity. The equicontinuity allows us to take a subsequence s.t. $\g_n$ converges uniformly to $\tilde{\g}$ as a capacity-parameterized curve (which is also given by Theorem~2 in  \cite{friz2015}). Together with the uniform convergence of $\l_n$ to $\l$, elementary calculus allows us to conclude that $\tilde{\g}$ is the Loewner transform of $\l$.
\end{proof}  

Before proving the proposition, let us review some further general background material on Loewner chains: The Loewner chain generated by a continuous function has local growth, and there are simple examples of continuous driving functions such that the Loewner chain is not a simple curve as well as examples where it is not even generated by a continuous curve. For instance, Lind, Marshall and Rohde (see \cite{lind2010collisions}) exhibit a driving function with  $\l (t) \sim 4\sqrt{1-t}$ when $t \to 1-$, 
that generates a Loewner chain with infinite spiral at time $1-$.  
An interesting subclass of simple curve is the quasiarc, which is the image of a segment under a quasiconformal mapping of $\m{C}$.  Marshall, Rohde \cite{marshall2005loewner} studied
when a continuous function generates a quasiarc (as a radial Loewner chain).
 Lind \cite{lind2005sharp} derived the following sharp condition for a driving function to generate a 
 quasislit half-plane in the chordal setting, \ie the image of $\m{H}\backslash [0,i]$ by a $K$-quasiconformal mapping fixing $\m{H}$ and $\infty$, its complement in $\m{H}$ is a quasiarc which is not tangent to the real line.

Let $T < \infty$, recall that $\Lip_{1/2}(0,T)$ is the space of H\"older continuous functions with exponent $1/2$, which consists of functions $\l(t)$ satisfying 
$$\abs{\l(s)- \l(t)} \leq c \abs{s-t}^{1/2}, \quad \forall 0\leq s,t \leq T$$
for some $c < \infty$, equipped with the norm $\norm{\l}_{1/2,[0,T]}$ being the smallest $c$.  

\begin{theorem}[\cite{lind2005sharp}] If the domain $\m{H} \backslash \g[0,T]$ generated by $\l$ is a quasislit half-plane,  then $\l \in \Lip_{1/2}(0,T)$.
Conversely, if $\l \in \Lip_{1/2}(0,T)$ with $ \norm{\l}_{1/2,[0,T]} < 4$, then $\l$ generates a $K$-quasislit half-plane for some $K$ depending only on $\norm{\l}_{1/2,[0,T]}$. 
\end{theorem}
The constant $4$ is sharp because of the spiral example mentioned before (where the local $1/2$-H\"older norm is as close to $4$ as one wishes).   
We remark that $H$ injects into $\Lip_{1/2}([0,\infty))$ since  
$$\abs{\l(t_1) - \l(t_2)} \leq \int_{t_1}^{t_2} \abs{\dot{\l}} \dd t \leq (t_2 -t_1)^{1/2} \left(\int_{t_1}^{t_2} \abs{\dot{\l}}^2 \dd t\right)^{1/2} \leq (t_2 -t_1)^{\frac{1}{2}} (2I(\l))^{1/2}.$$

The proof in \cite {lind2005sharp} is based on the following lemma, which gives a necessary and sufficient condition for being a quasislit half-plane in terms of the conformal welding, and is similar to Lemma~2.2 in \cite{marshall2005loewner}. 
More precisely, let $\l$ be the driving function of $\g$ with $ \norm{\l}_{1/2,[0,T]} < 4$, and let $f_T$ be the conformal mapping of $\m{H} \backslash \g[0,T] \rar \m{H}$, such that $f_T(z)-z$ is bounded, and $f_T(\g_T) = 0$. 
It is shown in \cite[Lemma~3]{lind2005sharp} that $f_T^{-1}$ extended to $\partial \m{H}$, and that for $0\leq t < T$, $\g_t$ has exactly two preimages of different sign which give a pairing on a finite interval of $\m{R}$. Then we extend the pairing to $f_T(x)$ and $f_T(-x)$.   The conformal welding $\phi : \m{R} \rar \m{R}$ is the decreasing function sending one point to its paired point.   

\begin{lemma}\label{lem_6_sharp}
    $\m{H} \backslash \g[0,T]$ is a $K$-quasislit half-plane if and only if $\phi$ is well-defined as above, and there exists $M > 0$ depending only on $K$, such that
    \begin{enumerate}[\upshape (i)]
       \item for all $x>0$,
       $$\frac{1}{M} \leq \frac{x}{-\phi(x)} \leq M,$$
       \item and for all $0 \leq x < y < z$ with $y-x= z-y$,
       $$\frac{1}{M} \leq \frac{\phi(x) - \phi(y)}{\phi(y) - \phi(z)} \leq M.$$
    \end{enumerate}
    And if $ \norm{\l}_{1/2,T} < 4$, the conformal welding of $\phi$ of $\g[0,T]$ satisfies both conditions, with $M$ depending only on $\norm{\l}_{1/2,[0,T]}$.
\end{lemma}

Quasiconformal mappings fixing $\m{H}$ can always be extended to a homeomorphism of the closure $\overline {\m{H}}$, thus induces a homeomorphism on $\partial \m{H}$. 
\begin{lemma}[\cite{lehto2012univalent} section 5.1] \label{lem_distorsion}
There exists a function $l(K)$ such that the boundary value $h$ of a $K$-quasiconformal mapping leaving $\m{H}$ and $\infty$ invariant is  always $l(K)$-quasisymmetric. \ie
   $$\frac{1}{l(K)} \leq \frac{h(x) - h(y)}{h(y) - h(z)} \leq l(K),$$
   for all $x<y<z$ and $y-x= z-y$. 
\end{lemma}

And conversely, 

\begin{theorem}[Kenig-Jerison extension, \cite{astala2008elliptic} Theorem~5.8.1] \label{lem_KJ}
There exists a function $K(k)$ such that 
every  $k$-quasisymmetric function $h$ on $\m{R}$ with  $h(0) = 0$ can be extended to a $K$-quasiconformal mapping with $i\m{R}_+$ fixed.  
\end{theorem}

The condition $h(0)=0$, and the extension fixes $i\m{R}_+$ is in particular convenient for the following proof. The extension exists for arbitrary quasisymmetric $h$ by translation. Now we are ready for proving Proposition~\ref{prop_finite_curve}.

~\\

\begin{figure}
\centering
\includegraphics[width=0.85\textwidth]{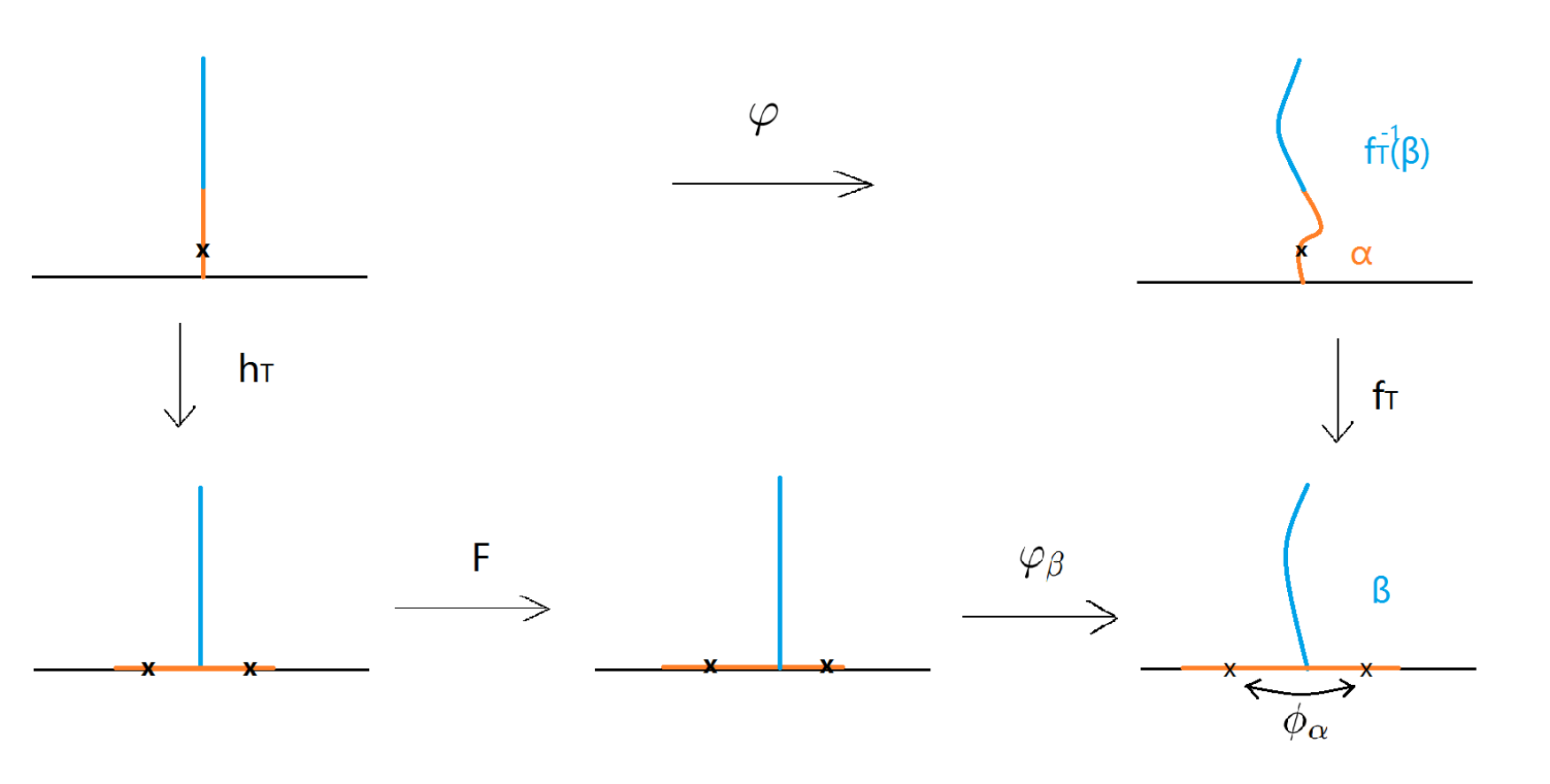}
\caption{\label{fig_concat} $h_T$ and $f_T$ are conformal mapping-out functions, $\varphi_\b$ is $K_\b$-quasiconformal compatible with $\b$, and $F$ is constructed in the proof. The mapping $\varphi$ which makes the diagram commute, then extended to $\m{H}$, is quasiconformal and compatible with $f_T^{-1}(\b) \cup \a$.} 
\end{figure}

\begin{proof} [Proof of Proposition~\ref{prop_finite_curve}] Let us first consider the case of driving functions $\l$ such that  $ L := I (\l) < 8$ (therefore $\norm{\l}_{1/2} < 4$). 
Thus for every $n \in \m{N}^*$, there exists $t_n >0$ and a $K(L)$-quasiconformal mapping $\varphi_n$, preserving $\m{H},0$ and $\infty$ such that $\g_{[0,n]} = \varphi_n (i[0,t_n])$, and $\varphi_n(i) = \g_1$. 
The module of the ring domain $\m{H}\backslash i[1,t_n]$ is bounded by $K(L)$ times the module of $\m{H} \backslash \g[1,n]$ since the latter is the image by the $K(L)$-quasiconformal map $\varphi_n$. The module of $\m{H} \backslash \g[1,n]$ converges to $0$, because both boundary components have spherical diameters bounded away from $0$ but has mutual spherical distance goes to $0$. For the proof see \cite{lehto1973quasiconformal}, p.~34. Thus $t_n \rar \infty$. 
The sequence $\varphi_n$ has a locally uniformly convergent subsequence, which does not contract to a boundary point because of $\varphi_n(i) = \g_1$. Let $\varphi$ denote the limit of this subsequence.   It is a $K(L)$-quasiconformal mapping, and since for every $y > 0$, $\varphi_n(iy)$ is on $\g$ for sufficiently large $n$ we can conclude that $\g = \varphi (i\m{R}_+)$.

The case of general driving functions $\l$ is treated by concatenating several pieces with small energy. This idea is used 
in \cite{lind2010collisions}, Thm~4.1. We repeat it here in order to see that the quasiconformal constant after concatenation is bounded by a constant which depends only on the constant of both parts of the quasiconformal curves. 

Let $\b$ be a $K_\b$-quasiconformal curve from $0$ to $\infty$ in $\m{H}$, $\varphi_\b$ a corresponding quasiconformal mapping. Let $\a := \g[0,T]$ a Loewner chain driven by $\l$ with 
$I_T (\l) \le 4$, with the centered mapping-out function $f_T$.  By Lemma~\ref{lem_6_sharp}, $\m{H}\backslash \a$ is $K_\a$-quasislit half-plane, where $K_\a$ depends only on $I_T(\a)$. It suffices to construct a quasiconformal mapping $\varphi$ compatible with  $f_T^{-1} (\b) \cup \a$,  with constant depending only on $K_\a$ and $K_\b$ (Figure~\ref{fig_concat}).

In fact, we only need to construct $F$ quasiconformal in $\m{H}$ which keeps track of the welding, \ie $$\varphi_{\beta} \circ F(-x) = \phi_\a \circ \varphi_\b \circ F(x),$$ 
where $\phi_\a$ is the conformal welding of $\a$  with the associated constant $M_\a$ as defined in Lemma~\ref{lem_6_sharp}.
The mapping $\varphi$  which makes the diagram commute can be extended by continuity to $\m{H} \to \m{H}$, and is the quasiconformal map that we are looking for. To construct $F$, first define $F$ on $\m{R}$ by   
$$F(x) = \varphi_\b^{-1} \circ \begin{dcases}
 \varphi_\b (x) & x\geq 0 \\ 
\phi_\a \circ \varphi_\b(-x) & x < 0,
\end{dcases}$$

Since $\varphi_\b^{-1}$ is a $K_\b$-quasiconformal mapping $\infty$ onto itself, it extends to  a $l(K_\b)$-quasisymmetric function on $\m{R}$ (Lemma~\ref{lem_distorsion}).    
 Moreover, the mapping
$$\psi: x \mapsto \begin{dcases}
 \varphi_\b (x) & x\geq 0 \\ 
\phi_\a \circ \varphi_\b(-x) & x < 0,
\end{dcases}$$
is $K$-quasisymmetric on each side of $0$, where $K$ depends only on $K_\a$ and $K_\b$. Besides, the ratio of dilatations on both sides is controlled by $M_\a$, \ie 
for $x > 0$, 
$$\frac{\psi(x)}{-\psi(-x)} = \frac{\varphi_\b (x)}{-\phi_\a \circ \varphi_\b(x)} \in \left[\frac{1}{M_\a}, M_\a \right]$$
thanks to Lemma~\ref{lem_6_sharp}. We conclude that $\psi$, hence $F$, are $K$-quasisymmetric with $K$ depending only on $K_\a$ and $K_\b$.  Using Theorem~{\ref{lem_KJ}}, we conclude that $F$ can be extended to a quasiconformal mapping preserving $i\m{R}_+$ with constant of quasiconformality depending only on $K_\a$ and $K_\b$. 

Cutting the driving function into small energy pieces, with $1/2$-Lip norm less than $4$, allows us to conclude the proof.
\end{proof}

\section{Reversibility via SLE$_{0+}$ large deviations}\label{sec_minimizing_curve_finite}

\subsection{SLE background}

We now very briefly review the definition and relevant properties of chordal and radial SLE.  The presentation is minimal and these properties are stated without proof.  
For further SLE background, readers can also refer to \cite{lawler2008conformally},
\cite{werner2004st_flour}.

  When $\k > 0$, the Loewner transform of  $(\sqrt{\k}B_t)_{t\geq 0}$ (where $B$ denotes a standard one-dimensional Brownian motion) is a random variable taking values in $\mc{L}_{\infty}$ and is called the \emph{chordal $\SLE_\k$ in $\m{H}$ from $0$ to $\infty$}.   
They are the unique processes having paths in $\mc{L}_{\infty}$ and satisfying \emph{scale-invariance} and the \emph{domain Markov property}. \ie for $u \in (0, \infty)$, the law is invariant under the scaling transformation
$$ (K_t)_{t\geq 0}  \mapsto (K^u_t  := u K_{u^{-2}t})_{t\geq 0}$$ 
 and for all $s\in[0,\infty)$, if one defines
$K_t^{(s)} =  K_{s, s+t} - \l_s$,
where $(K_s)$ is driven by $(\l_s)$, 
then the process $(K_t^{(s)})_{t\geq 0}$ has the same distribution as $(K_t)_{t \geq 0}$ and is independent of $\mc{F}_s = \s(\l_r: r\leq s)$. 

The scale-invariance makes it possible to define SLE$_\k$ also in other simply connected domains up to linear time-reparametrization (just take the image of the SLE in the upper half-plane via some conformal map from 
$({\m H}, 0, \infty )$ onto $(D, a, b)$ to define SLE from $a$ to $b$ in $D$. 

We are going to use the following known features of Schramm Loewner Evolutions: 

\begin{theorem}[\cite{Schramm2005basic}]\label{thm_transient_sle}
For $\k\in [0,4]$, $\SLE_\k$ is almost surely a simple curve $(\g_t)_{t \geq 0}$ starting at $0$.  Almost surely $\abs{\g_t} \rar \infty $ as $t \rar \infty$. 
\end{theorem}

\begin{theorem} [\cite{zhan2008}]
   For $\k \in [0,4]$, the distribution of the trace of $\SLE_\k$ in $\m{H}$ coincides with its image under $z \rar -1/z$.  
\end{theorem}

When $\l \in C_0([0,\infty))$ generates a simple curve $\g$ from $0$ to $\infty$ in $\m{H}$ (which is the case for SLE$_\kappa$ curves when $\kappa \le 4$),  then $\m{H} \setminus \g$ has two connected components $H^+ (\g)$ and $H^- (\g)$ having respectively $1$ and $-1$ on the boundary, that we can loosely speaking call the right-hand side and the left-hand side of $\g$. 
We also say that $\g$ is \emph{to the right (resp. left) of a subset $K\subset \m{H}$}, if $K \subset H^-(\g)$ (resp. $K \subset H^+(\g)$).

\begin {theorem} [\cite{Perco_Sch}] 
\label {Schrammevent}
For a point $z = x + iy \in {\m H}$, then for $w = x / y$  and $\kappa \le 4$, we define $h_\kappa(w)$ to be the probability that the SLE$_\kappa$ trace passes to the right of $\{z\}$. Then, 
$$h_\k(w) = \frac{1}{2}\int_w^{\infty} (s^2+1)^{-4/\k} \dd s / \int_0^{\infty} (s^2+1)^{-4/\k} \dd s.$$
\end {theorem}

 In fact, we will only use the fact that all these theorems hold for $\kappa$ very small. 

\subsection{One-point constraint large deviations and the minimizing curve}\label{sec_one_point}

For the record, let us now recall the standard large deviation theorem for Brownian paths. 
Let $T \in \m{R}_+$ and $C_0([0,T])$ be the set of continuous function $\l$ on $[0,T]$  with $\l(0) = 0$,  endowed with the $L^{\infty}$ norm, and $I_T$ as defined in Sec.~\ref{sec_Loewner_energy}. Then: 
\begin{theorem}[Schilder](See eg. \cite{Zeitouni}) 
The law of the sample path of $\sqrt{\k}B_t$ (the scaled Wiener measure $W_\k$) satisfies the large deviations principle with good rate function $I_T$ while $\k$ approaches $0$. 
More precisely, any closed set $F$ and any open set $O$ of $C_0([0,T])$,
$$\varlimsup_{\k \rar 0} \k \ln W_\k ( F ) \leq -\inf_{\l\in F} I_T(\l),$$
$$\varliminf_{\k \rar 0} \k \ln W_\k ( O ) \geq -\inf_{\l \in O} I_T(\l);$$
and that $I_T$ is lower semi-continuous and $\{\l \in  C_0([0,T]) \st I_T(\l) \leq c\}$ is compact for every $c \geq 0$.
\end{theorem}

Schilder's theorem is also valid when $T= \infty$ with rate function $I$,  but for $C_0([0,\infty))$ endowed with the norm $\sup_{t>0} |\l(t)|/(t+1)$ (see eg. \cite{Deuschel1989large}). In this paper we only use it for $T<\infty$.

We now fix a point $z$ in the upper half-plane with argument $\t$, and define 
$$D(z) := \{\l \in C_0([0,\infty)) \st \tau(z) < \infty \text{ or } \lim_{t \to \infty} \arg (f_t (z))  = \pi \},$$
where $(f_t)_{t\leq 0}$ is the centered flow driven by $\l$ (here and in the sequel, we  always choose the argument of a point in the upper half-plane to be in $(0, \pi)$). 

By Lemma~3 in \cite{Perco_Sch}, we see that $\l \in D(z)$ if and only if $z \notin H^+(\g)$.  Hence the probability that SLE$_\k$ is in $D(z)$ is $h_\k(w)$,  where $w = Re(z)/Im(z)$ with the notation in Theorem \ref{Schrammevent}. We also consider the function $F_\k(\cdot):=h_\k (\cot(\cdot))$.
We study this probability in the small $\k$ limit. 

\begin{prop}\label{prop_gd}
If $\t \in (0, \pi /2)$, then as $\k \to 0$, $-\k \ln F_\k(\t)$ converges to the infimum of $I( \l)$ over $D(z)$.  Furthermore,
this infimum is equal to $ - 8\ln(\sin(\t))$ and there exists a unique function $\l$ in $D(z)$ with this minimal energy.
\end{prop}

\begin{proof} The first fact and the existence of minimizing curve follow from the large deviation principle. It  is in fact a subcase of the multiple point constraint result that we will prove in Prop.~\ref{prop_gd_multiple}, so we do not give the detailed argument here.  In order to evaluate the value of this limit, we can use Theorem \ref{Schrammevent}:  
Indeed, for $w = \cot(\t) \geq 0$,
$$
	\k \ln(h_\k(w))  = \k \int_0^w \frac{h_\k'(s)}{h_\k(s)}\dd s + \k \ln (h_\k(0)),
	$$
and it is straightforward (see appendix) to check that if one defines $F(w) = 8w/(w^2+1)$, then the quantity 
$$\eps_\k (w) := -\k \frac{h'_\k(w)}{h_\k(w)} -F(w) = \frac{\k(w^2 +1)^{-4/\k}}{\int_w^{\infty}  (s^2+1)^{-4/\k} \dd s} - F(w) $$
converges to $0$ as $\k \to 0$, uniformly on $[0,\infty) $. Hence, 
 $$ \lim_{\k \to 0}   \k \ln(h_\k(w)) = -8 \int_0^w s ds / (1 + s^2) = - 4 \ln (w^2 + 1) = 8 \ln (\sin (\t)).$$

 Let us now prove the uniqueness of the minimizer (we will in fact also construct it explicitly):
  Suppose that $\g$ is a minimizing curve, denote the driving function of $\g$ by $\l$, with centered flow $(f_t)$. Let $z_t = f_t(z)$ and $w_t = \cot(\arg(z_t))$. 
 It is easy to see that $\tau (z)$ is finite. Indeed, if this is not the case, then $w_t \to - \infty$ when $t \to \infty$ and by the continuity of $t \mapsto w_t$, there exists $\tau < \infty$ s.t. $ w_{\tau} =0$ and $z_{\tau} \in i\m{R}_+$. But, $\l$ must be constant after $\tau$ to maintain the minimal energy. It would imply that $\gamma$ hits $z$, which contradicts the contrapositive assumption. 
 Note that for $t < T$, the curve $f_t(\g[t,\infty))$ is also a minimizing curve corresponding to $D(z_t)$  and  $\l$ is constant after the time $\tau(z)$. Hence, by the previously derived result, we get that for all $t <  \tau (z)$, $$\frac{1}{2}\int_t^{\tau(z)} \dot{\l}(s)^2 ds =  - 8\ln(\sin(\arg (z_t))).$$  
  By differentiation with respect to $t$,  we get that
\begin{equation}\label{eq_diff_f_dot}
 {(\dot{\l_t})^2} = 16 \times \frac { \Re (z_t)} { \Im (z_t) }  \times  \partial_t (\arg (z_t)) 
\end{equation}
But by Loewner's equation:
$$\partial_t \arg (z_t)= \Im ( (\partial_t z_t) / z_t )  = \Im ( 2/z_t^2  - \dot{\l}_t  / z_t ) $$
and a straightforward computation then shows that (\ref{eq_diff_f_dot}) can be rewritten as 
\begin{equation}\label{eq_f_dot}
 ( \dot{\l}_t - ( 8 \Re (z_t)   /  |z_t|^2  ))^2 = 0.
\end{equation}
This shows that if there exists a minimizer, then it is unique, as the previous equation describes uniquely $\lambda$ up to the hitting time of $z$. 

Conversely, if we define the driving function $\lambda$ that solves this differential equation for all times before the (potentially infinite) hitting time $\tau$ of $z$, we have indeed a minimizer:
One way to see that this $\l$ generates a curve passing through $z$ is to consider its the image in the unit disk $\m{D}$ by applying the conformal mapping $\psi_\t$ sending $z$ to $0$, $0$ to $e^{i\t}$, hence sending $\infty$ to $e^{-i\t}$. 
After the change of domain, the equation (\ref{eq_f_dot}) gives a simple characterization in the unit disc: 
The image of the minimizing curve $\g$ under $\psi_\t$ is symmetric with respect to the real axis. The part above the real axis can be viewed as a radial Loewner chain (in the radial time-parameterization) starting from $e^{i\t}$ driven by the function $\xi(t)$  such that $\cos(\xi(t)) = e^{-t}\cos(\t)$, hence hits $0$ in the $t \to \infty$ limit. The corresponding time in the original half-plane parameterization is finite, and the energy of $\g$ is indeed $-8\ln(\sin(\t))$.
\end {proof}

\begin{figure}
\centering
\includegraphics[width=0.38\textwidth]{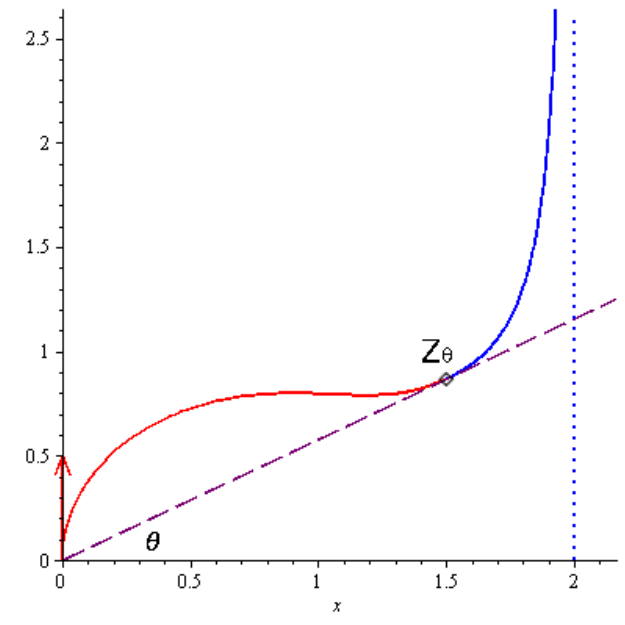}
\includegraphics[width=0.4\textwidth]{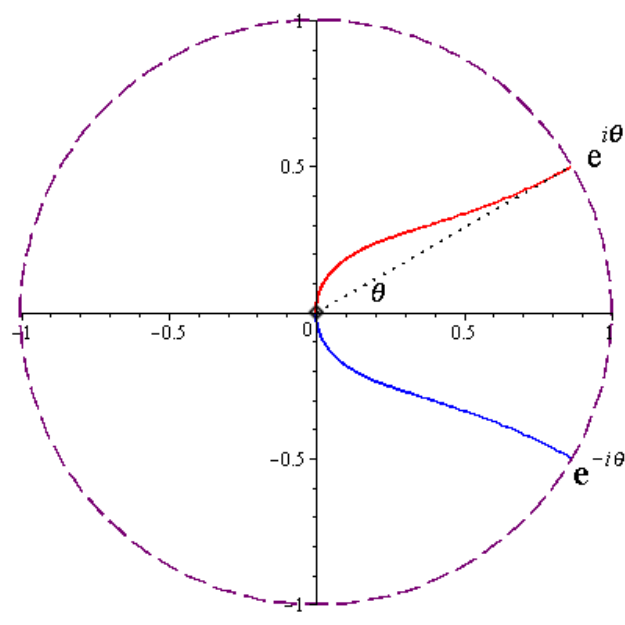}
\caption{\label{fig_curve_D} Sketch of the minimizing curve in $\m{H}$ (passing to the right of $z=  \rho e^{ i \pi /6}$ and of the
optimal curve in $\m{D}$ from $e^{-i\pi /6}$ to $e^{i\pi /6}$ (passing to the left of $0$). } 
\end{figure}

We list some remarks on the minimizer and its relation to SLE:
\begin{itemize}

\item
Readers familiar with the SLE$_\k(\rho)$ may want to interpret this first part of this curve in the radial characterization (from $e^{i\t}$ until it hits the origin) as a radial
SLE$_0(2)$ curve starting at $e^{i\t}$ with marked point at $e^{-i\t}$. By coordinate change (\cite{Schramm2005sle_change}) when coming back to the half-plane setting, 
the minimizing curve before hitting $z$ can be viewed as a chordal SLE$_0(-8)$ starting at the origin, with inner marked point at $z$. 

\item
The uniqueness of the minimizer $\g$ also implies  that $\g$ is the limit of SLE$_\k$ conditioned to be in $D(z)$, as $\k \to 0$. 
More precisely, if $W_\k^z$ denotes the conditional law of SLE$_\k$ in $D(z)$, then for all positive $\delta$, the $W_\k^z$-probability that the curve 
is at $\mc{L}_\infty$  distance greater than $\d$ of $\g$ goes to $0$ as $\kappa \to 0$.

\item
Note finally that Freidlin-Wentzell theory also provides a direct proof of the convergence of conditional SLE$_\k$ to the minimizer $\g$ on any compact time interval before the hitting time $T$ of $z$ by $\g$. Indeed, under the conditional law $W_\k^z$, the driving function $X_t$ and the flow of $z$ for $t < T$ are described by
 \begin{equation*}
(E_\k): \quad
\begin{dcases}
\ dX_t = \sqrt{\k} d\b_t + \frac{F(w_t)}{\Im (z_t) }dt + \frac{\vare_\k(w_t)}{\Im (z_t) } dt \\
\ dz_t = \frac{2}{z_t} dt - dX_t 
\end{dcases} 
\end{equation*}
where $\b$ is a Brownian motion and
$$
    w_t = \Re (z_t) / \Im (z_t),  \ F(w) =\frac{8w}{w^2+1} \wedge 0  \text{ and } \vare_\k(w) 
  \xrightarrow[\k \rar 0]{\text{unif. w.r. to }w \in \m{R}} 0
$$
as in the proof of Proposition~\ref{prop_gd}. The minimizer $\g[0,T]$ is driven by the solution to the deterministic differential equation $E_0$.  We can show that there exists a unique strong solution to $E_\k$ up to time $T$ with initial conditions $X_0 = 0, z_0 = z$. That solution to $E_\k$ converges to the unperturbed one in probability as $\k \to 0$, similar to the Freidlin-Wentzell theorem (see \cite{freidlin2012random}).

\end{itemize}

\subsection{Finite point constraints}\label{sec_finite_point}

Now we deal with the minimizing curve under finitely many point constraints. 
Let $\ad{z}$ denote a set of $n$ labeled points $\ad{z} = \{(z_1, \eps_1), \cdots, (z_n, \eps_n)\}$ in $ \m{H} \times \{ -1, 1 \}$. The label $\vare_i = +1 $ (resp. $\vare_i = -1$) is interpreted as ``right'' and ``left''.

Similarly to the one point constraint case, we define the set of functions \emph{compatible} with $\ad{z}$ as
$$D(\ad{z}) =  \{\l \in C_0([0,\infty)) \st \forall i \leq n,  \text {either }\tau_i < \infty \text{ or } \lim_{t \to \infty} \vare_i w_i(t) =  +\infty\},$$
where $\tau_i = \tau(z_i)$, $w_i(t) = \cot(\arg(f_t(z_i)))$ for $t<\tau_i$.

In the case where $\l$ generates a simple curve  $\g$ from $0$ to $\infty$, then $\l \in D(\ad{z})$ if and only if for every $i$, the point $z_i$ is not in $H^{-\eps_i} (\g)$. 

Let $D^t(\ad{z}) \subset C_0([0,t])$ denote the set of functions such that $\l(\cdot \wedge t) \in D(\ad{z})$. We also identify $\lambda \in C_0([0,t])$ with $\lambda( \cdot \wedge t) \in C_0([0,\infty))$ in order to make sense of hitting times, note that the comparison between $\tau_i$ and $t$ only depends on the restriction of the function to $[0,t]$, thus does not depend how is the function extended to $\m{R}_+$. Note also that a function of $C_0([0,t])$ is in $D^t ( \ad {z})$ if and only if for all $i \le n$, 
either $\tau_i \le t$ or $w_i (t) \vare_i \geq 0$. 

The main goal of this subsection is to prove the following result,  that 
states that the probability of SLE$_\k$ is in $D(\ad {z})$ decays exponentially as $\kappa \to 0$ with rate that is equal to the infimum of the energy in $D(\ad{z})$:

\begin{prop}\label{prop_gd_multiple}
There exists a positive constant $\tau$ such that for all $t \ge \tau$,
\begin{equation}
   \lim_{\k \rar 0} -\k \ln(W_\k (D(\ad{z}))) = \inf_{\l \in D(\ad{z})} I (\l) = \inf_{\l \in D^t(\ad{z})}  I_{t} (\l).
\end{equation}
Moreover, there exists a function $\l$ in $D^{\tau}(\ad{z})$ such that $I(\l)$ is equal to this infimum of the energy over $D(\ad{z})$.
\end{prop}

\begin{figure}
\centering
\includegraphics[width=0.95\textwidth]{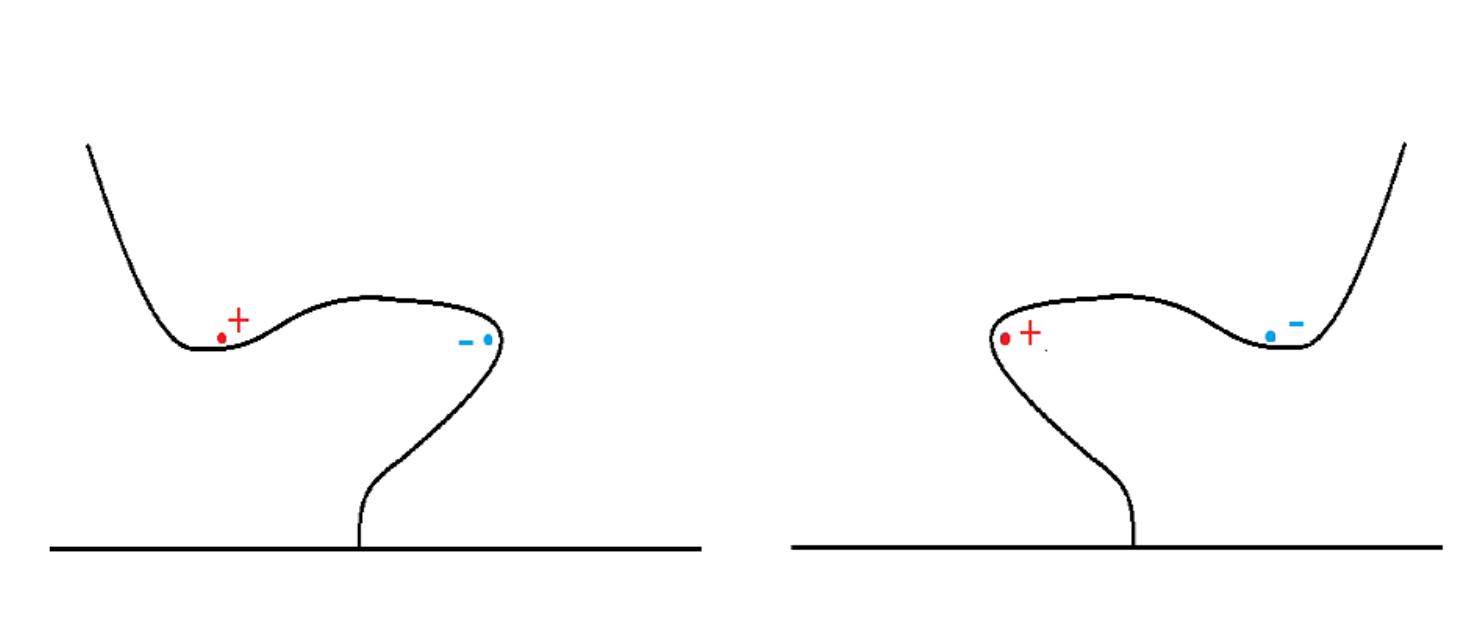}
\caption{\label{fig_symmetry} two curves compatible with the two-point constraint.} 
\end{figure}

Note that in this multiple point case, the minimizing curve is not necessarily unique anymore. For instance, we can consider two different points $z_1$ and $z_2$ that are symmetric to each other with respect to the imaginary axis, the left one assigned with $\vare_1= (+1)$ and the right one assigned $\vare_2 = (-1)$. Then, if there exists a minimizing curve in $D(\ad{z})$, the symmetric one with respect to 
the imaginary axis is necessarily a different curve, and is also minimizing the energy in $D(\ad{z})$ (see Figure~\ref{fig_symmetry}).

It is useful to consider the sets $D^{T}(\ad{z})$, as it  allows to reduce the study of the driving function $\l$ to the truncated function $\l\vert_{[0,T]}\in C_0([0,T])$, for which one can use large deviations easily. Moreover, it also shows that it suffices to look at the truncated part of the driving function on a finite interval to decide whether it is a minimizer of $D(\ad{z})$. 

By scale-invariance of the energy, it is enough to consider that case where all the points $z_i$ are in  $\m{D} \cap \m{H}$. 
For all $R > 0$, the hitting time $\tau_R= \inf\{t \geq 0 \st \abs{\g_t} = R\}$ is bounded. As half-plane capacity is non-decreasing, it follows that
$$2\tau_R = {\rm hcap}(\g[0,\tau_R]) \leq {\rm hcap}(R\m{D} \cap \m{H}) = R^2  {\rm hcap}(\m{D} \cap \m{H}).$$
This last half-plane capacity is in fact easily shown to be equal to $1$, so that $\tau_R$ is anyway bounded by $R^2/2$.

We have seen that if $\g$ has finite energy, then it is a simple curve. If $\tau_i > \tau_R$, then 
by comparing the harmonic measure seen from $z_i$ at both sides of $\g[0,\tau_R]$, we get that $\arg (f_{\tau_R}(z_i))  \in (0,\t') \cup (\pi-\t',\pi)$, for some  $\t'= \t' (R)$ that tends to $0$ as
$R \to \infty$. But the energy needed after $\tau_R$ to change the sign of $w_i(t)$ will therefore be larger than $-8\ln(\sin(\t'))$ (which goes to infinity as $R \to \infty$).

But on the other hand, Lemma \ref {lem_finite_energy} will tell us that one has an \emph{a priori} upper bound on the infimum of the energy in $D(\ad{z})$.
Hence, we can conclude that there exists $R$ such that
for any function in $D(\ad{z})$ with energy smaller that $ \inf_{\l \in D(\ad{z})} I(\l) +1$, the sign of $w_i$ can not change after $R^2 /2$. This therefore implies that $\l ( \cdot \wedge R^2/2)$ is also in 
$D ( \ad {z})$.
 
\begin{lem} \label{lem_finite_energy}
There exists a simple curve $\g$ with finite-energy that visits all $n$ points $z_1, \ldots, z_n$.
\end{lem}

\begin{proof}
   In the one-point constraint case, we have seen that the energy needed for touching a point in $\m{H}$ is finite. It is then easy to iteratively use this to define a curve with finite energy that touches successively all constraint points (start with the part of the minimizing curve that hits $z_1$, and then continue in the complement of this first bit with the minimizing curve that hits $z_i$ where $z_i$ has the smallest index among points that have not been hit yet and so on). We choose the driving function to be constant after this curve has hit all the $n$ points.  
\end{proof}

Let us now define the set $O^T(\ad{z})$ of driving functions $\l$ in $C_0 ([0,T])$ such that for all $i \le n$, 
$ T < \tau_i$ and the sign of the real part of $f_T ( z_i)$ is $\vare_i$. We can note that this set is open in $C_0([0,T])$. 
Indeed, if $\l$ is in $O^T ( \ad {z})$, then by inspecting the evolution of the $n$ points under the Loewner flow, if we perturb only slightly the driving function (in the sense of the $L^\infty$ metric on $[0,T]$), we  will not change the fact that one stays in $O^T ( \ad {z})$. Note of course that $O^T ( \ad {z}) $ is a subset of $D^T ( \ad {z})$. The complement of $D^T(\ad{z})$ is the union of $O^T(z_i,-\vare_i)$, so that $D^T (\ad{z})$ is closed in $C_0([0,T])$.

The following lemma will be useful in our proof of Proposition \ref {prop_gd_multiple}:

\begin{lem}\label{lem_equal}  For every $t > 0$,
  $$\inf_{\l \in D^t(\ad{z})}  I_{t} (\l) = \inf_{\l \in  {O^t}(\ad{z})}  I_{t} (\l).$$
\end{lem}

In the proof of the lemma, we will do surgeries on driving functions of the following type: 
saying that the part $\l[t,t+\d]$ is replaced by $a[0,\d]$ where $a(0)=0$, means that the new driving function $\tilde{\l}$ is defined as, $\tilde{\l}(s) = \l(s)$ for  $s\in [0,t]$, $\tilde{\l}(t+s) - \tilde{\l}(t)= a(s)$ for $s\in [0, \d]$, and $\tilde{\l}(s') - \tilde{\l}(s)= \l(s') - \l(s)$ for $s,s' \geq t+\d$.  

Our goal is to show that for any $\l$ in $D^t ( \ad {z})$ with finite energy $L$, and any $c> 0$, we can find  a perturbed $\l_{\vare} \in {O^t}(\ad{z})$ with energy 
$I (\l_{\vare}) \leq L +c$, and $\norm{\l- \l_{\vare}}_{\infty} \leq c$.

Note that (just as in the previous argument showing that $O^t ( \ad {z})$ is open), if $\g$ does not go through any point $z_i$, then one can ensure by taking $\vare$ small enough that any such 
perturbation will not change on which side $z_i$ ends. Hence, it suffices to prove the result in the case where $\g$ in fact visits all the points $z_1, \ldots, z_n$.

   The idea is that when the flow starting from $z_i$ is about to hit $0$ (ie. when $\g$ is about to hit $z_i$), one can replace a small portion of the driving function by some optimal curve (targeting at a well chosen point) to avoid a neighborhood of $z_i$, up to the time when the point $z_i$ tends to the real line but away from $0$, and therefore becomes very hard to be reached again. The modification on energy is controlled as well as the impact on other constraints point.  
\\

\begin{proof}[Proof of Lemma~\ref{lem_equal}]

  Let us only explain how to proceed in the case with two constraint points $z_1, z_2$ and $\eps_1 = \eps_2 = -1$, as the other 
  cases are treated similarly and require no extra idea. 
  Assume that $\g$ visits first $z_1$ and then $z_2$ (at times $T_1$ and $T_2$) and that $\eps_1= \eps_2 = -1$.

     For the simplicity of notation we write $I_{t,T} (\g) := I_T(\g) - I_{t}(\g)$. We use the notation $f(t)$ for the centered flow  instead of $f_t$ to avoid too heavy subscripts.    As before, let $z(t) = f(t)(z_1)$ and $w(t):= \cot(\arg(z(t)))$. \\

 Using a scaling we assume $\abs{f(T_1)(z_2)} > 1$.  For every $0 < \eta \ll 1$, there exists $ \d < \eta$, such that $I_{T_1-s, T_1+s}(\g) \leq \eta$, $\abs{f(T_1-s)(z_1)} \leq \eta $ and $\abs{f(T_1 \pm s)(z_2)} \geq 1$ for all $s \leq \d$.
 
 For every $s \leq \d$,  we know from Section \ref{sec_one_point} that the minimal energy for a curve to hit a point with argument $\t$ is $-8\ln(\sin(\t)) = 4 \ln(w^2 +1)$ where $w = \cot(\t)$. 
   Hence $w(T_1 -s)^2 \leq \exp(\eta/4) - 1$.  For every $\vare > 0$, we choose a point $\tilde{z}_1$ with argument slightly smaller than $\arg(z(T_1 -s))$, and close enough to $z(T_1-s)$, such that the minimizing curve driven by $a$, targeting at $\tilde{z}_1$ is to the right of  a small neighborhood of $z(T_1-s)$, and such that $I(a) \leq  \eta + \vare$.
   
   Let $s_1$ denote the hitting time of $\tilde{z}_1$ under the driving function $a$, then $s_1 \asymp s$.  
   If we do the surgery $a[0, \infty)$ to the initial driving function from time $T_1 -s$, the cotangent of the argument $\tilde{w}$ driven by the new function, satisfies that $\tilde{w}(T_1 -s + s_1) < 0$ by construction. After $T_1 -s + s_1$, the driving function is constant, so $\tilde{w}(t)$ shrinks very fast towards $-\infty$ under the flow driven by the $0$ function which is $z \mapsto \sqrt{z^2+4t}$. 
   When $\abs{\tilde{w}(t)}$ reached the threshold of ``the point need energy $2L$ to be hit by the rest of the curve'',  \ie when $\tilde{w}(t)^2 \geq \exp(L/2) -1$, let $\tilde{T_1}$ denote that time, then for any surgery made after time $\tilde{T_1}$, $\forall r \geq \tilde{T_1}$, $\tilde{w}(r) < 0$ if the energy $I_{\tilde{T_1}, \infty}$ is smaller than $2L$.  Thus we don't worry about the sign of $\tilde{w}(t)$ to change in the future if the energy on the remaining curve is not too much perturbed. Moreover the time $\tilde{T_1} \to T_1$ as $s\to 0$. 
   
   Now choose $s < \d$, such that $\tilde{T_1}\leq  T_1 + \d$. We do the surgery of $\l$ by replacing the part $\l([T_1-s, \tilde{T_1}])$ by $a([0, \tilde{T_1} - T_1 + s ])$. The energy of the total curve is increased by at most $\vare$.
   
   The new curve $\g^1$ driven by $\l^1$ does not necessarily pass through $z_2$. Write the new flow $(f^1(t))$, since we assumed $\abs{f(T_1-s)(z_2)} \geq 1$, by Gronwall type argument, the impact of replacing a piece of length $ s_1 + s_2 $ of driving function, will change $f(\tilde{T_1})(z_2)$ to $\tilde{z}_2 :=f^1(\tilde{T_1})(z_2)$, with Euclidean distance less than 
   $$\int_{0}^{ s_1 + s_2 } \abs{\dot{\l}(T_1-s + t) - \dot{a}(t)} \dd t 
   \leq 
   \sqrt{ s_1 + s_2 } \sqrt{2 I_{T_1-s, \tilde{T_1}} (\l)+ 2 I (a) }
   \leq 
   \sqrt{\d} \sqrt{4 \eta +2\vare },$$
   which can be chosen to be arbitrarily small. 
   
   Now we apply the same argument to $z_2$. 
   Since $\l^1$ and $\l$ have same increments on $[\tilde{T_1},\infty)$, $f^1(\tilde{T_1}) (\g^1[\tilde{T_1},\infty))$ passes through $f(\tilde{T_1})(z_2)$ at time $T_2 - \tilde{T_1}$. With arbitrairily small compromise in the energy, we can modify again the $\l^1$ to $\l^2$ such that the curve $f(\tilde{T_1})(\g^1[\tilde{T_1},\infty))$ passes to the right of a neighborhood of $f(\tilde{T_1})(z_2)$, which contains $\tilde{z}_2$ by well choosing $\eta, \d, \vare$ in the previous step. 
\end{proof}

~\\

We are now finally ready for proving Proposition~\ref{prop_gd_multiple}. Let $T$ be a function $C_0([0,\infty)) \to \m{R}_+$, we define the set of functions $\mc{D}^{T}(\ad{z})$  \emph{compatible with $\ad{z}$ up to time $T$} 
 to be $\{\l \in C_0([0,\infty)), \l|_{[0,T(\l)]} \in D^{T(\lambda)}(\ad{z})\}$. 
We keep the notations at the beginning of Section~\ref{sec_finite_point} and those before Lemma~\ref{lem_finite_energy}.
~\\

\begin{proof}[Proof of Proposition~\ref{prop_gd_multiple}] Let $M < \infty$ be the energy of the curve $\g$ passing through all constraint points constructed in Lemma~\ref{lem_finite_energy}, $T(\g)$ be the hitting time of the last point of $\g$. We choose $\t'$ s.t.  $-8\ln(\sin(\t')) \geq  M+1$. There exists $R \in \m{R}_+$, such that for every $i \leq n$, if $\tau_i > \tau_R$, $f_{\tau_R}(z_i)$ has argument in $(0,\t')\cup (\pi-\t', \pi)$ and let $\tau = R^2/2$. By construction, $T(\g)<\tau_R(\g) \leq  \tau$. So $D^{\tau}(\ad{z})$ is non empty, and the infimum of the energy of curves in $D^{\tau}(\ad{z})$ is less than $M$.

We know that SLE$_\k$ is almost surely transient and does not hit $\ad{z}$. Thus apart from a set $N$, with zero $W_\k$-measure for all $\k \leq 4$, the symmetric difference between $D(\ad{z})$ and $\mc{D}^{\tau_R}(\ad{z})$ consists of driving functions which give different signs to $w_i(\tau_R)$ and the limit of $w_i(t)$ when $t\rar \infty$ for some $i \leq n$, 
\begin{equation*}
D(\ad{z}) \D \mc{D}^{\tau_R}(\ad{z})  =  \{\l \in C_0([0, \infty)) \st \exists i \leq  n, \tau_i = \infty \text{ and } w_i(\tau_R)w_i(t) \xrightarrow[t\rar \infty]{} -\infty\} \cup N
\end{equation*}
where $N$ is $W_k-$null set for all $\k \le 4$.

By the domain Markov property, $w_i(\tau_R)w_i(t) \rar -\infty$ means that the Loewner curve reaches a point with argument $\t'$ or $\pi-\t'$ starting from time $\tau_R$. Hence the probability of $\sqrt{\k}B$ stays in the symmetric difference $D(\ad{z})\Delta \mc{D}^{\tau_R}(\ad{z})$ is bounded by $2 n F_\k(\t')$, where $F_\k(\t)$ is the probability that SLE$_\k$ is to the right of $z$ with argument $\t$. 
 The inequality $\tau_R \leq \tau$ holds for all driving functions, gives the upper-bound of the probability on the symmetric difference between $D(\ad{z})$ and $\mc{D}^{\tau}(\ad{z})$:
\begin{align*}
W_\k (D(\ad{z}) \D \mc{D}^{\tau}(\ad{z})) \leq & W_\k (D(\ad{z}) \D \mc{D}^{\tau_R}(\ad{z})) \\
& + W_\k (\exists i \leq  n \st w_i(\tau_R)w_i(t) \rar +\infty, w_i(\tau)w_i(t) \rar -\infty) \\
\leq  &  W_\k (D(\ad{z}) \D \mc{D}^{\tau_R}(\ad{z}))+W_\k (\exists i \leq  n \st w_i(\tau_R)w_i(t) \rar -\infty) \\
 \leq &2W_\k (D(\ad{z}) \D \mc{D}^{\tau_R}(\ad{z})) \leq  4 n F_\k(\t').
\end{align*}

 The inequality also holds for $t \geq \tau$.
Using Schilder's theorem on $[0, t]$ and the fact that $D^t(\ad{z})$ is non empty and closed and  $O^t(\ad{z}) \subset D^t(\ad{z})$ is open,  we have for any $\vare>0$, 
\begin{align*}
\k \ln W_\k(D(\ad{z})) \leq& \k \ln \left[W_\k( \mc{D}^t(\ad{z})) + 4nF_\k(\t') \right]  \\
&=  \k \ln  \left[W_\k( D^t(\ad{z})) + 4nF_\k(\t') \right] 
\leq -\inf_{\l\in  D^t(\ad{z})} I_{t}(\l) + \vare\\
\k \ln W_\k ( D(\ad{z})) \geq& \k \ln \left[W_\k({O^t}(\ad{z}))  - 4nF_\k(\t')\right]
\geq -\inf_{ \l \in  {O^t}(\ad{z})} I_{t}(\l) -\vare,
\end{align*}
for $\k$ is sufficiently small, since $\lim_{\k \to 0}\ln ( F_{\k}(\t')) = M+1$ and $\inf_{\l\in  D^t(\ad{z})} I_{t}(\l) \leq M$.
 By Lemma~\ref{lem_equal}, we can therefore conclude that for all $t \ge \tau$,
  $$\lim_{\k \rar 0}-\k \ln W_\k ( D(\ad{z}))  =\inf_{\l \in  D^t(\ad{z})} I_{t}(\l). $$
The set $D^{\tau}(\ad{z})$ is closed and $I_{\tau}^{-1}([0,M])$ is compact and non empty,
so that there exists $\l_0 \in D^{\tau}(\ad{z})$ such that $I(\l_0) = \inf_{t\geq \tau} \inf_{\l \in  D^{t}(\ad{z})} I_{t}(\l) =  \inf_{\l \in  D(\ad{z})} I(\l)$.
\end{proof}

\subsection{Proof of reversibility of the energy}
We can now conclude the proof of Theorem~\ref{main1}. 
For a curve $\g$ from $0$ to $\infty$ driven by $\l$, the reversed curve $-1/\g$ is driven by a continuous function by  Theorem~\ref{thm_Loewner_Kufarev}, we define it as $\Rev(\l)$.  In particular, by Proposition~\ref{prop_finite_curve}, the functional $\Rev$ is well-defined for finite energy functions.  Let $-1/\ad{z}$ be the constraint set obtained by taking the image of $\ad z$ by $z \to -1/z$  equipped with the opposite assignment.  Write $Z(\g)$ for the collection of all finite constraint sets compatible with $\g$. We have in particular $\ad{z} \in Z(\g)$ if and only if $-1/\ad{z} \in Z(-1/\g)$.

Let ${\mc I} ( \ad {z})$ denote the infimum of  $I ( \l)$ over 
$\l \in D(\ad{z})$ ie. it is  the minimal energy needed to fulfill the constraint $\ad{z}$.

\begin{prop} \label{prop_rev}
Let $\g_0$ be a simple curve from $0$ to $\infty$ in $\m{H}$. The energy of $\g_0$ equals to the supremum of $\mc{I}(\ad{z})$ over  $\ad{z}\in Z(\g_0)$. 
\end{prop}

This proposition implies Theorem~\ref{main1}. Indeed, by reversibility of SLE$_\k$, we know that
$W_\k (D(\ad{z})) = W_\k(D(-1/\ad{z}))$. Thus by Proposition~\ref{prop_gd_multiple}, ${\mathcal I} ( \ad{z} )  =
 {\mathcal I} (  -1/ \ad{z})$. Hence (writing $\ad {z'} = {-1/\ad {z}}$), 
 $$ I ( -1 / \g_0 ) =  \sup_{\ad{z}' \in Z(-1/\g_0)} {\mathcal I} ( \ad{z}' )  =
 \sup_{\ad{z} \in Z(\g_0)} {\mathcal I} (  -1/ \ad{z}) 
 = \sup_{\ad{z} \in Z(\g_0)} \mathcal {I} ( \ad{z}) 
 = I (\gamma_0) .$$

\begin{proof}[Proof of Proposition~\ref{prop_rev}]
  Let $L(\g_0) =\sup_{\ad{z}  \in Z( \g_0)} {\mathcal I} (  \ad{z})$. As $L(\g_0) \leq  I(\g_0)$, we only have to prove that $L(\g_0) \geq I(\g_0)$ in the case where $L( \g_0)$ is finite.
 
 Let $(\ad{z_n})_{n \geq 1}$ be a increasing sequence of finite point constraints, with assignments compatible with $\g_0$ and  such that $\bigcup_{n\geq 1} \{ z_n \} = \m{Q} \times \m{Q}_+$. We pick a minimizer $\g_n$ 
 of the energy in $D ( \ad{z_n}) $ arbitrarily (the existence of such minimizers follows from Proposition~\ref{prop_gd_multiple}). In addition, we know that $I(\g_n)$ is non-decreasing in $n$, and bounded from above
 by $L(\g_0)$.  
In view of Corollary~\ref{cor_compact} (with the same notation as in that corollary), there exists a subsequence $n_k$ such that $\varphi_{n_k}$ converges uniformly locally to some $\varphi$.
Note that the energy of  $\tilde{\g} := \varphi(i \m{R}_+)$ can not be larger than $L(\g_0)$.  
 
 We  now show that $\tilde{\g}$ is compatible with $\ad{z_n}$ for all $ n$.  Assume for instance that the point constraint on $z_n$ is $+1$. Then, for $n_k \ge n$,  $\varphi_{n_k}^{-1}(z_n)$ has non-negative real part. Again by quasiconformality of $\varphi_{n_k}^{-1}$,  we know that a subsequence of $(\varphi_{n_k}^{-1})_{k \geq 1}$ has local uniform limit $\varphi^{-1}$ so that the real part of  $\varphi^{-1}(z_n) $ is also non-negative. Thus it follows that $\tilde{\g}$ is compatible with all the constraints $\ad{z_n}$.
Given that both $\g_0$ and $\tilde{\g}$ are simple curves, it follows readily that they are equal, which in turns implies that   $I ( \gamma_0) =  I(\tilde{\g}) \le L ( \gamma_0)$. 
\end{proof}

\section {Comments}
\label {commentssection}

In the present final section, we will provide two more direct derivations of results that will provide some insight and information about the energy of a Loewner chain. In fact, it might look at 
first sight that each of these results could be used in order to get an alternative derivation of Theorem \ref {main1}, but it does not seem to be the case for reasons that we will comment on. 

\subsection{Conformal restriction}\label{sec_conformal_restriction}

 In this subsection we study the variation of the energy of a given curve when the domain varies (this is similar to the conformal restriction-type ideas initiated in \cite {lawler2003conformal} in the SLE-framework,  similar computations can be also found in the Section 9.3 of \cite{Dub{\'e}dat2007commutation}).
Let $K$ be a compact hull at positive distance to $0$.  The simply connected domain $H_K : = \m{H} \backslash K$ coincides with $\m{H}$ in the neighborhoods of $0$ and $\infty$.
Let $\g \subset \m{H}$ be a finite energy curve connecting $0$ and $\infty$, driven by $W$ with centered flow $(f_t)$, and assumed to be at positive distance to $K$. 

\begin{prop}\label{prop_energy_change_finite}
The energy of $\g$ in $(\m{H},0 ,\infty)$ and in $(H_K, 0 , \infty)$ differ by 
$$
   I_{H_K, 0 , \infty}(\g) - I(\g) = 3\ln\abs{\psi_0'(0)} + 12 m^l(\g, K; \m{H}), 
$$
where $\psi_0$ is the conformal equivalence $H_K\rar \m{H}$ fixing $0,\infty$ and $\psi_0(z)$ being normalized as $z +O(1)$ near $\infty$; $m^l(A, B; \m{H})$ is the measure of Brownian loops in $\m{H}$ intersecting both $A$ and $B$ for $A,B \subset \m{H}$.
\end{prop}

The Brownian loop measure is a natural conformally invariant measure on the unrooted loops defined in \cite{lawler2004}. In particular we have $ m^l(\g, K; \m{H}) =  m^l(-1/\g, -1/K; \m{H})$. 

We also remark that 
$\abs{\psi_0'(0)}$ is the probability of a Brownian excursion in $\m{H}$ starting from $0$ not hitting $K$ (this was first observed in \cite{Virag2003beads}, see also \cite{werner2004st_flour} Lemma~5.4; 
by a Brownian excursion in $\m{H}$ we mean a process in $\m{H}$ which is Brownian motion in $x$-direction and an independent $3$-dimensional Bessel process in $y$-direction, see \cite {wernerconformal} for more information). 
Note that if $\g$ is a analytic curve such that there exists some $K$ where $\g$ is the hyperbolic geodesic in $(H_K, 0, \infty)$, then $I_{H_K, 0, \infty}(\g)= 0$. Hence, 
$$I(\g) = -3\ln\abs{\psi_0'(0)} - 12 m^l(\g, K; \m{H}),$$
which does provide an interpretation of the energy $I(\g)$ as a function  of ``conformal distance'' between $K$ and $\g$, and shows directly reversibility of the energy for such curves $\g$.
However, such analytic curves form a rather small class of curves (for instance, the beginning of the curve fully determines the rest of it) so that it does not seem to  allow to deduce the general reversibility result easily from this fact.
~\\

Proposition~\ref{prop_energy_change_finite} is a consequence of the following proposition as $t \to \infty$.

\begin{prop}  Let $K_t := f_t(K)$, and $\psi_t: \m{H} \backslash K_t \rar \m{H}$ the conformal mapping fixing $0$ and $\infty$, such that $\psi_t(z)-z$ bounded.  
   $$I(\psi_0(\g[0,t])) =  \frac{1}{2} \int_0^t \left[\dot{W}_s - \frac{3 \psi_{t}''(0)}{\psi_{t}'(0)}\right]^2 \dd t. $$
 Moreover,
   $$I(\psi_0(\g[0,t])) - I(\g[0,t]) = 3\ln\abs{\psi_0'(0)} + 12 m^l(\g[0,t], K; \m{H}) - 3\ln\abs{\psi_t'(0)},$$  
   with $3\ln\abs{\psi_t'(0)} \rar 0$ when $t \rar \infty$. 
\end{prop}
The first equality is obtained by computing the driving function of $\psi_0(\g)$ as in \cite{lawler2003conformal} Section~5 and then the difference is identified using the identity: 
$$m^l(\g[0,t], K; \m{H})  = -\frac{1}{3} \int_0^t S\psi_s (0) \dd s,$$
where $S$ is the Schwarzian derivative. The identification is explained in \cite{lawler2003conformal} Section~7.1. Moreover, $-S\psi_s (0)/6$ can be interpreted as the half-plane capacity of $K_s$ seen from $0$ as we explained in Section~\ref{sec_Loewner_energy}. 
~\\

\subsection{Two curves growing towards each other} \label{sec_two_curves}

The commutation relations for SLE \cite{Dub{\'e}dat2007commutation} are closely related to the coupling of both ends of the SLE curves used to show SLE reversibility.
In the present subsection, we will derive directly  an analogous commutation relation for the energy of two curves growing towards each other. One may again wonder whether it is possible to 
 deduce the reversibility of the Loewner energy from this commutation, but it seems again to be strictly weaker than Theorem~\ref{main1}.  

\begin{figure}
\centering
\includegraphics[width=0.95\textwidth]{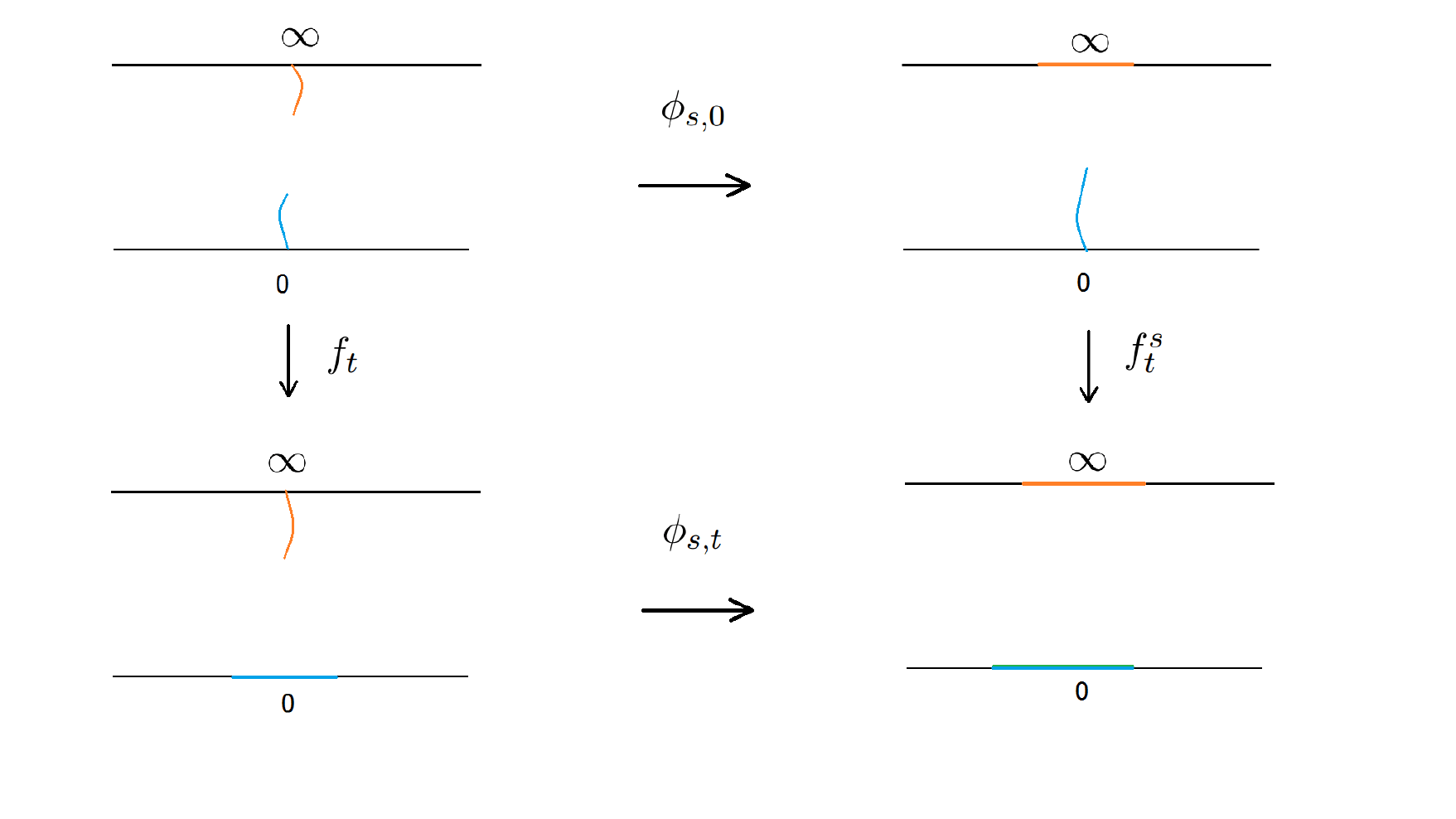}
\caption{\label{fig_two_curves} $f_t$ and $f_t^s$ are centered mapping-out functions normalized at $\infty$, and $\phi_{s,0}$ is the centered mapping-out function normalized at $0$. The mapping $\phi_{s,t}$ makes the diagram commute.} 
\end{figure}

 Let us compute the energy of a curve $\g^T$ in $\m{H}\setminus \tilde{\g}^S$, where $\g^T:=\g[0,T]$ starting at $0$ is driven by $W$ and $\tilde{\g}^S:= \tilde{\g}[0,S]$ is driven by $U$ in $(\m{H}, \infty, 0)$ at positive distance of $\g^T$, \ie the image of $(\tilde{\g}_s)$ by $z\mapsto-1/z$ is driven by $U$. The centered flow $(\phi_{s,0})_s$ associated with $\tilde{\g}^S$ is, for every $s \leq S$, the mapping-out function of $\tilde{\g}^s$ sending the tip to $\infty$, fixing $0$  and being normalized as $\phi_{s,0}'(0) = 1$.
In particular $(\phi_{s,0})_s$ satisfies the Loewner equation:
$$\partial_s(\phi_{s,0}) (z) = -\phi_{s,0}(z)^2 (2\phi_{s,0}(z) + \dot{U_s}), \quad \phi_{0,0}(z) = z.$$
Thus 
$$\partial_s (\phi_{s,0}'(z)) =  -6\phi_{s,0}(z)^2 \phi_{s,0}'(z)-2 \phi_{s,0}(z) \phi_{s,0}'(z) \dot{U_s}$$
\begin{align*}
\partial_s (\phi_{s,0}''(z)) = &-12 \phi_{s,0}(z)(\phi_{s,0}'(z))^2- 6\phi_{s,0}(z)^2\phi_{s,0}''(z) \\
&-2(\phi_{s,0}'(z))^2 \dot{U_s}-2\phi_{s,0}(z)\phi_{s,0}''(z)\dot{U_s}
\end{align*}
Since $\phi_{s,0}(0) = 0$ and $\phi_{s,0}'(0) = 1$,
\begin{align*}
\partial_s \left(\frac{\phi_{s,0}''(0)}{\phi_{s,0}'(0)}\right) &= \frac{\partial_s\phi_{s,0}''(0)}{\phi_{s,0}'(0)} - \frac{\phi_{s,0}''(0)\partial_s\phi_{s,0}'(0)}{(\phi_{s,0}'(0))^2} =-2 \dot{U_s} 
\end{align*}
which implies
$\phi_{s,0}''(0)/\phi_{s,0}'(0) =  -2 U_s$. 
The map $\phi_{s,t}$ is defined as illustrated in the Figure~\ref{fig_two_curves}, where $f_t$ and $f_t^s$ are the (time reparametrized) centered flow of $\g^T$ and $\phi_{S,0}(\g^T)$. Since $\phi_{s,t}''(0)/\phi_{s,t}'(0)$ does not depend on the normalization at the image,  we deduce that 
$$\phi_{s,t}''(0)/\phi_{s,t}'(0) =  -2 U_s^t$$
where $(U^t_s)_{s\geq 0}$ is the driving function for $f_t(\tilde{\g})$ parameterized by capacity of $\tilde{\g}$.

\begin{lem} With above notations,
   \[
  I_{\m{H} \setminus \tilde{\g}^S, 0, \tilde{\g}_S}(\g[0,T]) = \frac{1}{2} \int_0^T \left[\dot{W}_t - \frac{3 \phi_{S,t}''(0)}{\phi_{S,t}'(0)}\right]^2 \dd t
 =\frac{1}{2} \int_0^T \left[\dot{W}_t + 6 U_S^t\right]^2 \dd t.
\]
    Furthermore,
    \begin{align*}
    I_{\m{H} \setminus \tilde{\g}^S, 0, \tilde{\g}_S}(\g^T) - I(\g^T) &= 12 m^l(\g^T, \tilde{\g}^S; \m{H}) - 3\ln (\phi'_{S,T}(0)) \\
    &= 12 m^l(\g^T, \tilde{\g}^S; \m{H}) - 3\ln
    \left(
    \frac{H(\m{H} \backslash \g^T \cup \tilde{\g}^S, \g_T, \tilde{\g}_S)H(\m{H},0,\infty)}
    {H(\m{H} \backslash \g^T, \g_T, \infty) H(\m{H} \backslash \tilde{\g}^S, 0,\tilde{\g}_S)}
    \right), \end{align*}
  where $H(D,x,y)$ is the Poisson excursion kernel in a domain $D$ at boundary points $x$ and $y$ with respect to analytic local coordinates in the neighborhood of each of them. We choose local coordinates to coincide if it is the same neighborhood involved. 
\end{lem}
The quotient term of Poisson excursion kernel does not depend on the local coordinates (for more discussions on the Poisson excursion kernel, readers may refer to \cite{Dub{\'e}dat2009sle} Section~3). 
The first equality is obtained as in the previous section, and the second one is due to the computation above. In particular, we see how the energy of $\tilde{\g}^S$ redistributes on $\g^T$ in $(\m{H} \setminus \tilde{\g}^S, 0, \tilde{\g}_S)$ in a rather intricate way. 

Also notice that the last expression of the difference shows the symmetry in $\g^T$ and $\tilde{\g}^S$. Hence we get the following commutation relation:
\begin{cor}
The sum of the energy of one slit and  the energy of the second slit in the domain left vacant by the first one, targeting at the tip of the first slit, does not depend on the ordering of the slits:
$$ I(\tilde{\g}^S) +I_{\m{H} \setminus \tilde{\g}^S, 0, \tilde{\g}_S}(\g^T)
 = I(\g^T) + I_{\m{H} \setminus \g^T, \infty, \g_T}(\tilde{\g}^S).$$
 \end{cor}

 To conclude, let us  remark that  the above corollary is also a simple corollary of Theorem~\ref{main1}: Let us define the hyperbolic geodesic $\g$ in the two-slit domain between the two tips.
 The concatentation of this geodesic with the two slits defines a curve from the origin to infinity in the upper half-plane. Then
 $$
 I(\g^T \cup \g \cup \tilde{\g}^S) = I(\g^T) + I_{\m{H} \setminus \g^T, \infty, \g_T}(\g \cup \tilde{\g}^S) 
 = I(\g^T) + I_{\m{H} \setminus \g^T, \infty, \g_T}(\tilde{\g}^S), 
$$
 (the first equality is due to the additivity of $I$ and Theorem~\ref{main1} applied in $(\m{H} \setminus \g^T, \infty, \g_T)$). Considering the completed curve $\g^T \cup \g \cup \tilde{\g}^S$ in $(\m{H}, \infty, 0)$, we get
 $$I(\g^T \cup \g \cup \tilde{\g}^S) =I(\tilde{\g}^S) + I_{\m{H} \setminus \tilde{\g}^S, 0, \tilde{\g}_S}(\g^T),$$
 which provides an alternative proof of the corollary.

\appendix

\section{One point estimate}

\begin{lem}\label{lem_unif}
Let $F(w) =\frac{8w}{w^2+1} \wedge 0$, and $\vare_\k(w) = -\k \frac{h'_\k(w)}{h_\k(w)} - F(w)$, then
$$\vare_\k(w) = -\k \frac{h'_\k(w)}{h_\k(w)} -F(w) = \frac{\k(w^2 +1)^{-\frac{4}{\k}}}{\int_w^{\infty}  (s^2+1)^{-\frac{4}{\k}} \dd s} - F(w) \xrightarrow[\k \rar 0]{\text{unif. in } \m{R}} 0.$$
\end{lem}

\begin{proof}
For $w>0$, we can write
\begin{equation*}
  \vare_\k(w) = \frac{\k(w^2 +1)^{-\frac{4}{\k}}}
  {\int_w^{\infty}  (s^2+1)^{-\frac{4}{\k}} \dd s} -\frac{8w}{w^2 +1} =  \frac{2}{(w^2+1) M_\k(w)}  - \frac{8w}{w^2 +1}
\end{equation*}
after a change of variable $t^2 = \frac{(w^2+1)s^2}{s^2+1}$ in the integral:
$$M_\k(w) = \frac{2}{\k} \int_{w}^{\sqrt{w^2+1}} (w^2 +1 -t^2)^{\frac{4}{\k} - \frac{3}{2}} \dd t.$$
It can be bounded by
\begin{align*}
M_\k(w) &\leq \frac{1}{\k w}  \int_{w}^{\sqrt{w^2+1}} 2t(w^2 +1 -t^2)^{\frac{4}{\k} - \frac{3}{2}} \dd t \\
& = \frac{1}{\k w} \frac{2\k}{8- \k} [-(w^2 + 1 - t^2)^{\frac{4}{\k} - \frac{1}{2}}]_w^{\sqrt{w^2+1}} \\
& = \frac{2}{w(8-\k)}
\end{align*}
For $1> \vare > 0$,
\begin{align*}
M_\k(w) &\geq \frac{1}{\k \sqrt{w^2 + \vare}}  \int_{w}^{\sqrt{w^2+\vare}} 2t(w^2 +1 -t^2)^{\frac{4}{\k} - \frac{3}{2}} \dd t \\
& = \frac{1}{\k \sqrt{w^2 + \vare}} \frac{2\k}{8- \k} [-(w^2 + 1 - t^2)^{\frac{4}{\k} - \frac{1}{2}}]_w^{\sqrt{w^2+\vare}} \\
& = \frac{2}{(8-\k)\sqrt{w^2 +\vare}}(1- (1-\vare)^{\frac{4}{\k} - \frac{1}{2}}).
\end{align*}
Then we obtain bounds on $\vare_\k$:
$$\vare_\k(w) \geq \frac{-\k w}{w^2+1} \geq -\frac{\k}{2},$$
where we used $\frac{w}{w^2 +1} \leq \frac{1}{2}$. And
\begin{align*}
  \vare_\k(w) &\leq \frac{1}{w^2 +1} \left(\frac{(8-\k)\sqrt{w^2 +\vare}}{ 1- (1-\vare)^{\frac{4}{\k} - \frac{1}{2}}}-8w\right) \\
  & \leq \frac{1}{w^2 +1} \abs{\frac{(8-\k)\sqrt{w^2 +\vare} -8w + 8w (1-\vare)^{\frac{4}{\k} - \frac{1}{2}} }{ 1- (1-\vare)^{\frac{4}{\k} - \frac{1}{2}}}}\\
  &\leq \frac{1}{w^2 +1} \frac{(8-\k)(\sqrt{w^2 +\vare }- w) + \k w + 8w (1-\vare)^{\frac{4}{\k} - \frac{1}{2}}}{1- (1-\vare)^{\frac{4}{\k} - \frac{1}{2}}} 
\end{align*}   
for every $\k$, we choose $\vare$ such that $(1-\vare)^{\frac{4}{\k} - \frac{1}{2}}= \k$. One can check that $\vare \sim \frac{\k \ln(\k)}{4}$, and notice that $\sqrt{w^2 +\vare} - w \leq \sqrt{\vare}$. The upper-bound of $\vare_\k(w)$ becomes
$$\vare_\k(w) \leq \frac{1}{w^2 +1} \frac{(8-\k)\sqrt{\vare } + \k w + 8w \k}{1- \k} \leq \frac{(8-\k)\sqrt{k\ln(k)}}{1-\k} + \frac{9\k}{1-\k}\frac{1}{2},$$
which does not depend on $w$ and converges  to $0$ with speed $\sqrt{\k\ln(\k)}$. 

For $w\leq 0$, $\vare_\k(w)$ can be bounded brutally: 

$$0 \leq \vare_\k(w) = \frac{\k(w^2 +1)^{-\frac{4}{\k}}}
  {\int_w^{\infty}  (s^2+1)^{-\frac{4}{\k}} \dd s} 
\leq 
\frac{\k}{\int_0^{\infty}  (s^2+1)^{-\frac{4}{\k}} \dd s} 
= \vare_\k(0). 
$$
Using a change of variable $ u = \ln(s^2+1)/\k$, one gets
$$
\vare_\k(0) = \frac{2}{\int_0^{\infty}e^{-4u} e^{u\k}/\sqrt{e^{uk}-1}\dd u}
\leq 
\frac{2}{\int_0^1 e^{-4u}/ \sqrt{e^{uk}-1} \dd u} 
\leq
\frac{2 \sqrt{\k}}{\int_0^1  e^{-4u}/\sqrt{2u}\dd u}
= c \sqrt{\k}.  $$
 The second inequality holds for $\k < c'$, where $c, c' >0$. We conclude that $\vare_\k$ converges uniformly to $0$ with speed at least $\sqrt{\k}$ on $(-\infty, 0]$.
\end{proof}


\bigskip
\footnotesize
\noindent\textit{Acknowledgments.}
 I would like to thank Wendelin~Werner for numerous inspiring discussions as well as his help during the preparation of the manuscript. I am also grateful to Steffen~Rohde for his comments on the first draft,  to Julien~Dub{\'e}dat for nicely pointing out the similar idea in \cite{Dub{\'e}dat2007commutation}, to Yves~Le~Jan for conversations on the topic and to the referees for very helpful comments on the first version of this paper. This research was partly supported by SNF  (grant no. 155922).

\end{document}